\documentclass[letterpaper, 10 pt, conference]{ieeeconf}  % Comment this line out
\IEEEoverridecommandlockouts                              % This command is only
\usepackage{graphicx}           % to include graphics (pdflatex
\usepackage{amsmath,amssymb}    % math packages, e.g. for align and pmatrix environments
\usepackage[usenames]{xcolor}   % you might consider using the dvipsnames option for more colors
\usepackage{tikz}               % graphics inside latex; very powerful
\usepackage[expansion=false]{microtype} %typographic enhancements
\usepackage{amsthm}
\usepackage{mathtools}
\usepackage{amsopn}
\usepackage{varwidth} 
\usepackage{xcolor,colortbl}
\usepackage{tabu}
\usepackage{comment}
\usepackage{enumerate}
\usepackage{array}
\usepackage{csquotes}
\usepackage{tikz}
\usepackage{tikz-qtree}
\usetikzlibrary{trees,positioning,shapes}
\usepackage{pgfplots}
\usepackage{algpseudocode}
\usepackage{nicefrac}
\usepackage{bbm}
\usepackage[rightcaption]{sidecap}
\usepackage{import}
\usepackage{textcomp}
\usepackage{slashed}
\usepackage[section]{placeins}
\usepackage{tcolorbox}
\definecolor{grey}{rgb}{0.57, 0.64, 0.69}

\theoremstyle{plain}
\newtheorem{theorem}{Theorem}
\newtheorem{lemma}[theorem]{Lemma}
\newtheorem{proposition}[theorem]{Proposition}

\theoremstyle{definition}
\newtheorem{definition}[theorem]{Definition}

\newtheorem{remark}[theorem]{Remark}

\newtheorem{problem}{Problem}

\newcommand{\bbC}{\mathbb{C}}

\newcommand{\bbN}{\mathbb{N}}

\newcommand{\bbR}{\mathbb{R}}

\newcommand{\bA}{\boldsymbol{A}}
\newcommand{\bB}{\boldsymbol{B}}
\newcommand{\bC}{\boldsymbol{C}}
\newcommand{\bD}{\boldsymbol{D}}
\newcommand{\bE}{\boldsymbol{E}}

\newcommand{\bG}{\boldsymbol{G}}
\newcommand{\bH}{\boldsymbol{H}}
\newcommand{\bI}{\boldsymbol{I}}
\newcommand{\bJ}{\boldsymbol{J}}

\newcommand{\bL}{\boldsymbol{L}}
\newcommand{\bM}{\boldsymbol{M}}

\newcommand{\bP}{\boldsymbol{P}}
\newcommand{\bQ}{\boldsymbol{Q}}
\newcommand{\bR}{\boldsymbol{R}}
\newcommand{\bS}{\boldsymbol{S}}
\newcommand{\bT}{\boldsymbol{T}}
\newcommand{\bU}{\boldsymbol{U}}
\newcommand{\bV}{\boldsymbol{V}}

\newcommand{\bZ}{\boldsymbol{Z}}

\newcommand{\bPi}{\boldsymbol{\Pi}}
\newcommand{\bPsi}{\boldsymbol{\Psi}}
\newcommand{\bSigma}{\boldsymbol{\Sigma}}

\newcommand{\diff}{\,\mathrm{d}}
\newcommand{\quot}[1]{\enquote{#1}}
\DeclareMathOperator{\im}{im}

\DeclareMathOperator{\subsEq}{\subseteq}

\title{\LARGE \bf
Convex Synthesis of Accelerated Gradient Algorithms for Optimization and Saddle Point Problems using Lyapunov functions
}

\author{Dennis Gramlich, Christian Ebenbauer and Carsten W. Scherer% <-this % stops a space
\thanks{Dennis Gramlich and Christian Ebenbauer are with the Institute for Systems Theory and Automatic Control,
        University of Stuttgart, Germany
        {\tt\small ce@ist.uni-stuttgart.de}}%
\thanks{Carsten W. Scherer is with the Institute of Mathematical Methods in Engineering, Numerical Analysis and Geometric Modeling,
        University of Stuttgart, Germany
        {\tt\small carsten.scherer@mathematik.uni-stuttgart.de}}%
}

\begin{document}

\maketitle
\thispagestyle{empty}
\pagestyle{empty}

%%%%%%%%%%%%%%%%%%%%%%%%%%%%%%%%%%%%%%%%%%%%%%%%%%%%%%%%%%%%%%%%%%%%%%%%%%%%%%%%
\begin{abstract}

This paper considers the problem of designing accelerated gradient-based algorithms
for optimization and saddle-point problems. The class of objective functions 
is defined by a generalized sector condition.
This class of functions contains strongly convex functions with Lipschitz gradients
but also non-convex functions, which allows not only
to address  optimization problems but also saddle-point problems. 
The proposed design procedure
relies on a suitable class of Lyapunov functions and 
on convex semi-definite programming. 
The proposed synthesis allows the design of algorithms that reach the performance of state-of-the-art accelerated gradient methods and beyond.
\end{abstract}

%%%%%%%%%%%%%%%%%%%%%%%%%%%%%%%%%%%%%%%%%%%%%%%%%%%%%%%%%%%%%%%%%%%%%%%%%%%%%%%%
\section{Introduction}

Gradient-based optimization algorithms are a standard tool in science and engineering. 
Many of these algorithms take the form of feedback interconnection between a time-discrete linear system and the gradient of the objective function.
In case of a convex objective function, the corresponding gradient satisfies a certain sector condition. Hence such a feedback configuration falls in the class of so called Lur'e systems \cite{lur1944theory}, which have been extensively studied in control theory. 
In recent years, results from Lur'e systems and techniques from robust control theory have been exploited to analyze convergence rates and robustness of known optimization algorithms and to design novel algorithms. Some of those new publications rely on IQCs (integral quadratic constraints) from robust control to generate convergence results.
%As an example, IQCs were applied in \cite{lessard2016analysis} and \cite{michalowsky2019robust} to analyze the performance of numerical algorithms and design optimized methods.
For example, IQCs were used in \cite{lessard2016analysis} to find upper bounds for the convergence rates of existing algorithms. This work was later extended to synthesis of algorithms in \cite{lessard2019direct}.
These IQC-based approaches gave rise to the development of the Triple Momentum Method \cite{van2017fastest}. This method has the fastest known upper convergence bound for strongly convex functions with Lipschitz gradients.
Other related work that analyzes optimization algorithms from a dynamical systems perspective is for example given in \cite{fazlyab2018analysis}
and \cite{michalowsky2014}, 
%and \cite{hu2017dissipativity}
where also Lyapunov function techniques and robust control theory are employed, or in \cite{wilson2016lyapunov}, where discrete-time algorithms are analyzed based on continuous-time counterparts. 
In addition, semi-definite programming formulations have been proposed in  \cite{drori2014performance} and \cite{taylor2017smooth} to analyze the convergence properties of first order optimization methods.
Further related results are discussed in the recent paper \cite{michalowsky2019robust}, where the design of robust algorithms for structured objective functions based on IQC theory is considered.

In this paper, we address \emph{convex design} (convex synthesis) of gradient-based algorithms for optimization and saddle point problems, where the class of objective functions is defined by a generalized sector condition. In particular, the contributions of this paper are as follows.
First we consider classes of functions that are more general than the classes of strongly convex functions usually considered in the literature.
%mentioned above. 
In particular, the classes under consideration also contain non-convex functions, which we utilize in our procedure to design algorithms capable of searching for saddle points instead of minima. For example, the ability to search for saddle points allows us to apply the design method to optimization problems with equality constraints.
Second, based on a rather general class of Lyapunov functions, we derive convex synthesis conditions for algorithm design in the form of linear matrix inequalities.
Specifically, we provide a non-conservative convexification in the sense that
the analysis matrix inequalities (when algorithm parameters are given) are feasible if and only if the synthesis matrix inequalities (when algorithm parameters are decision variables) are feasible, i.e. our design procedure is not more conservative than the corresponding analysis.
This is in contrast to many other results in the literature, 
where the step from convex analysis to convex synthesis 
is only possible by imposing additional assumptions (such as
fixed IQC multipliers or quadratic Lyapunov functions).
 In the case of strongly convex functions, our design procedure reaches the same convergence rates as the Triple Momentum Method and it
 allows to incorporate additional structural properties of the objective function to 
 design tailored algorithms with even faster convergence rates, as demonstrated in the paper.

\section{Problem statement and preliminary results}

\subsection{Notation}

%We write $C^m(\bbR^d)$ for the set of $m$-times continuously differentiable functions from $\bbR^d$ to $\bbR$ and denote the gradient of a function $f \in C^1(\bbR^d)$ by $\nabla f$. 
%Matrices are denoted with bold, capital letters, while vectors are written in lower case letters.
By $\|v\|$, we denote the Euclidean norm of a vector $v \in \mathbb{R}^n$ and by $\|\bA\|$ the spectral norm of a matrix $\bA \in \mathbb{R}^{n \times n}$. The spectrum of a matrix will be denoted by $\sigma(\bA)$ and for the spectral radius we will write $\rho(\bA)$. We will also often use the notation $\|v\|_{\bA}^2 = v^T\bA v$ 
for the semi-norm defined by a positive semi-definite matrix $\bA$, which is a full norm whenever $\bA$ is positive definite. If $\bA_1,\bA_2$ are two symmetric matrices of the same dimensions, then we write $\bA_1 \succ \bA_2$ ($\succeq$) if $\bA_1 - \bA_2$ is positive (semi-) definite and $\bA_1 \prec \bA_2$ ($\preceq$) if $\bA_1 - \bA_2$ is negative (semi-) definite. With $\bA^\dagger$, we will denote the Moore-Penrose pseudo inverse of a matrix, while $\bA^T$ will denote its transpose. The orthogonal projection matrix onto the kernel of a matrix $\bA$ will be denoted by $\bPi_{\ker \bA}$. In large matrix equations, we will sometimes write $\bA^T\bB(\star)$. In that case, $(\star)$ is to be understood as a copy of the matrix $\bA$.
%Gradient based algorithms are interpreted as discrete linear systems, defined by matrices $\bA \in \bbR^{n\times n},\bB \in \bbR^{n \times d},\bC \in \bbR^{d\times n}$, in interconnection with the gradient of an objective function $f \in S \subsEq C^1(\bbR^d)$. Here, the natural number $n$ is always the dimension of the state space of the optimizer and the natural number $d$ is always the dimension of the domain of $f$. States of the optimizer are always written with an $x$, while elements from the domain of $f$ are always denoted by $z$.

\subsection{Problem statement}

Consider the gradient based algorithm defined by
\begin{align}
\label{eq:algo}
x_{k+1} &= \bA x_k + \bB \nabla f(\bC x_k),
\end{align}
where $x_k \in \bbR^{n}$ and the matrices $\bA \in \bbR^{n\times n},\bB \in \bbR^{n \times d},\bC \in \bbR^{d\times n}$ are the algorithm parameters to be designed.
The objective function $f \in C^1(\bbR^d)$ is assumed to satisfy the following generalized sector condition for all $z_1,z_2 \in \bbR^d$:
%The objective function is denoted by $f$ and we assume that $f$ satisfies the following generalized sector condition for all $z_1,z_2 \in \bbR^d$:
\begin{align}
	%\begin{split}
		\frac{1}{2}\| z_1 - z_2\|_{\bM}^2 &\leq f(z_2) - f(z_1) + (\nabla f(z_1))^T(z_1-z_2) \nonumber\\
		&\leq \frac{1}{2}\| z_1 - z_2 \|_{\bL}^2,
	%\end{split} 
	\label{eq:4SCLcharacterisation}
\end{align}
where $\bM \preceq \bL \in \bbR^{d\times d}$ are given symmetric matrices. 
In the following, $S(\bM,\bL)$ denotes the set of all $C^1$ functions that satisfy \eqref{eq:4SCLcharacterisation}.
Note that $S(m \bI_d, l \bI_d)$, $m<l$, is a set of strongly convex functions,
as typically found in the literature.
In the case $f \in C^2(\bbR^d)$, \eqref{eq:4SCLcharacterisation} is equivalent to $\bM \preceq H_f(z) \preceq \bL$ for all $z \in \bbR^d$, where $H_f$ denotes the Hessian of $f$.

The algorithm design problem addressed in this paper is formally stated as:

\begin{problem}
	\label{problem1}
    For given $n \ge d$, $\bM\preceq \bL$, and convergence rate $\rho \in [0,1[$,  we aim to \emph{design} matrices $( \bA, \bB,  \bC) \in  \bbR^{n\times n} \times \bbR^{n\times d} \times \bbR^{d\times n} $ such that for any $f \in S(\bM,\bL)$ there exist $x_f^* \in \bbR^n$ and $c_f \in \bbR_{\geq 0}$ such that
    \begin{align*}
        \nabla f(z_f^*) = 0 \mathrm{~for~} z_f^* := \bC x_f^*
    \end{align*}
    and the iterates $x_k$ of \eqref{eq:algo} satisfy
    \begin{align*}
        \| x_f^* -  x_k\| &\leq c_f \rho^k \| x_f^*- x_0\|,
    \end{align*}
    for any $x_0 \in \bbR^{n}$, $k \in \bbN_0$.
    %and such that there exists  
    %\begin{align*}
    %    \| x_f^* -  x_k\| &\leq c_f \rho^k \| x_f^*- x_0\|,
    %\end{align*}
    %for any $x_0 \in \bbR^{n}$, $k \in \bbN_0$. (We say then: \quot{The algorithm has convergence rate $\rho$}.)
\end{problem}

 In our setting, \emph{design} (synthesis) refers to computing the algorithm parameters $(\bA,\bB,\bC)$ by solving a  \emph{convex} optimization problem, i.e. a semi-definite program.
 
 Our goal is solving Problem \ref{problem1}. The following Problem \ref{problem2}  is similar to Problem \ref{problem1} with the slight modification that all the functions $f$ under consideration have their critical points in $z_f^* = 0$. This is favourable for the application of tools from robust control theory, which are often formulated for fixed-points in zero.

\begin{problem}
	\label{problem2}
	For given $n \geq d$, symmetric matrices $\widetilde{\bL} \succeq 0$ and $\bM$, and $\rho \in [0,1[$, design matrices 
	$(\widetilde \bA, \widetilde \bB, \widetilde \bC) \in  \bbR^{n\times n} \times \bbR^{n\times d} \times \bbR^{d\times n} $ satisfying the constraint
	\begin{align}
	\widetilde \bC(\widetilde{\bA} - \bI_n)^{-1} \widetilde \bB \bM = \bI_d \label{eq:Constraint}
	\end{align}
	such that for any $f \in S_0(0,\widetilde{\bL}) := \{ f\in S(0,\widetilde{\bL}): \nabla f(0) = 0 \}$ there exists $c_f \in \bbR_{\geq 0}$ such that the iterates of \eqref{eq:algo} satisfy
	\begin{align*}
	\| x_k\| &\leq c_f \rho^k \| x_0\|
	\end{align*}
	for any $x_0 \in \bbR^{n}$ and $k \in \bbN_0$.
\end{problem}

The subsequent theorem states that the two problems are equivalent.

\begin{theorem}
	\label{thm:ConvTrInv}
	Let symmetric matrices $\bM \preceq \bL$ be given, set $\widetilde{\bL}:= \bL - \bM$ and fix $\rho \in [0,1[$.
	Then the matrices $(\bA,\bB,\bC)$ solve Problem \ref{problem1} if and only if the matrices $(\widetilde{\bA},\widetilde \bB,\widetilde \bC)$ solve Problem \ref{problem2}, where $\widetilde{\bA} = \bA + \bB\bM\bC$, $\widetilde{\bB} = \bB$, $\widetilde{\bC} = \bC$.
\end{theorem}

This theorem justifies that we can solve Problem \ref{problem2} instead of Problem \ref{problem1}.

\subsection{Properties of the class $S(\bM,\bL)$}

%The function class $S(\bM,\bL)$ is a generalization of the class of strongly convex functions with Lipschitz continuous gradients $S(m \bI_d,l \bI_d)$, which is often considered in the literature.
%and which is defined as the set of functions $f \in C^1(\bbR^d)$ satisfying
%\begin{align*}
%\frac{m}{2}\|z_1 - z_2\|^2 &\leq f(z_2)-f(z_1) - (\nabla f(z_1))^T(z_2-z_1)\\ & \leq \frac{l}{2}\|z_1 - z_2\|^2.
%\end{align*}
%In the generalized case of $S(\bM,\bL)$, the scalar values $m,l$ with $0<m\leq l$ are replaced by symmetric matrices $\bM,\bL$ with $\bM \preceq \bL$. 
This subsection serves the purpose of introducing some important properties of $S(\bM,\bL)$.
The first result gives some equivalent characterizations for when $f \in S(0,\bL)$ holds true. Note, that these conditions can be applied to any class $S(\bM,\bL)$ by using the fact $f \in S(\bM,\bL) \Leftrightarrow (z\mapsto f(z) - \frac{1}{2}z^T\bM z) \in S(0,\bL-\bM)$.

\begin{lemma}[Characterizations for $f \in S(0,\bL)$]
	\label{lem:LipschitzConv}
	Let $\bL \succeq 0$ and $f \in C^1(\bbR^d)$. All conditions below, holding for all $z_1,z_2 \in \bbR^n$, are equivalent to $f\in S(0,\bL)$:
	\begin{enumerate}
	    \setlength\itemsep{1em}
		\item $0 \leq f(z_2)-f(z_1) - (\nabla f(z_1))^T(z_2-z_1) \leq \frac{1}{2}\|z_1-z_2\|_{\bL}^2$,
		\item $0\leq(\nabla f(z_1) - \nabla f(z_2))^T(z_1-z_2)\leq \|z_1-z_2\|_{\bL}^2$,
		\item $\frac{1}{2}\|\nabla f(z_1) - \nabla f(z_2)\|_{\bL^\dagger}^2
		\leq f(z_2) - f(z_1) + (\nabla f(z_1))^T(z_1 - z_2)$ and $\bPi_{\ker \bL} (\nabla f(z_1) - \nabla f(z_2)) = 0$,
		\item $\|\nabla f(z_1) - \nabla f(z_2)\|_{\bL^\dagger}^2
		\leq (\nabla f(z_1) - \nabla f(z_2))^T (z_1-z_2)$ and $\bPi_{\ker \bL} (\nabla f(z_1) - \nabla f(z_2)) = 0$.
	\end{enumerate}
\end{lemma}

Not all possible variations of matrices $\bM\preceq \bL$ should be considered for optimization. For example, if there exists a singular matrix $\bQ$ such that $\bM \preceq \bQ \preceq \bL$, then the function $f$ defined by $f(z) = \frac{1}{2}z^T\bQ z + v^T z$, where $v$ is not in the range of $\bQ$, would be an element of $S(\bM,\bL)$ without any critical point. Therefore this set $S(\bM,\bL)$ would not make sense as a set of objective functions, since we cannot solve Problem \ref{problem1} for it. The following Lemma characterizes when such cases can be avoided.

\begin{lemma}[Well-posed pairs $\bM, \bL$]
	\label{lem:4.1EigCond}
	Let $\bM,\bL \in \bbR^{d\times d}$ be symmetric matrices with $\bM \preceq \bL$. Then the following five statements are equivalent:
	\begin{enumerate}
		\item The matrices $\bM$ and $\bL$ have the same numbers of positive and negative, and no zero eigenvalues.
		\item Any symmetric matrix $\bQ \in \bbR^{n\times n}$ with $\bM \preceq \bQ \preceq \bL$ is non-singular.
		\item $\bL + \bM$ is non-singular and the spectral radius of $(\bL + \bM)^{-1}(\bL - \bM)$ is smaller than one.
		\item $\bM$ is non-singular and $\bM^{-1} \bL$ has only positive eigenvalues.
		\item $\bM$ and $\bL$ are non-singular and congruent, i.e. there exists a non-singular matrix $\bT \in \bbR^{d\times d}$ with $\bM = \bT^T\bL\bT$.
	\end{enumerate}
\end{lemma}

\begin{remark}
    In Lemma \ref{lem:4.1EigCond}, statement 1) serves the purpose of giving the reader a good intuition for the property under consideration. Statement 2) and 3) will be useful in later proofs. Note that in particular 2) prevents the counter-example we constructed in the motivation of this lemma. Statement 4) offers the most efficiently verifiable test of the considered property, by the fact that the verification whether a matrix has positive eigenvalues can be done by solving a Lyapunov equation.
\end{remark}

Because of the importance of this property we define a new notation for matrices $\bM,\bL$ fulfilling one and thus all conditions in Lemma \ref{lem:4.1EigCond}.

\begin{definition}[Loewner-congruence ordering on symmetric matrices]
	For symmetric matrices $\bM,\bL\in \bbR^{d\times d}$, we introduce the partial ordering
	\begin{align*}
	\bL \succeq_c \bM : \Leftrightarrow \begin{cases}
	\bL - \bM \text{ is positive semi-definite}\\
	\bL \text{ and } \bM \text{ are congruent}
	\end{cases}.
	\end{align*}
\end{definition}

Under the Loewner-congruence ordering, a critical point exists, is unique,
and a simple gradient method converges to the critical point, as stated in
the following results.

\begin{proposition}[A simple gradient method]
	\label{prop:4.1GradContraction}
	Let $\bL \succeq_c \bM$ be non-singular. Then for any convergence rate $\rho > \rho \left( (\bL + \bM)^{-1}(\bL - \bM) \right)$ there exists $r\in \bbR_{>0}$ such that
	\begin{align}
	z \mapsto z - 2(\bM + \bL)^{-1}\nabla f(z) \label{eq:StructuredGD}
	\end{align}
	is a contraction for all $f \in S(\bM,\bL)$ with contraction constant $\rho$ on the Banach space $(\bbR^d,\|\cdot\|_{\bP})$, where $\bP = (\bL + \bM)((\bL - \bM)^{\dagger} + r\bPi_{\ker (\bL-\bM)})(\bL + \bM)$.
\end{proposition}

\begin{remark}
    \label{rem:StructuredGD}
    For $\bM \preceq_c \bL$, the optimizer defined by $(\bA,\bB,\bC)$ with $\bA = \bC = \bI_d$ and $\bB = -2(\bL+\bM)^{-1}$ realizes the contraction in Proposition \ref{prop:4.1GradContraction}. As a consequence of the Banach fixed-point theorem, it converges faster than any convergence rate $\rho > \rho_\mathrm{grad}:= \rho \left( (\bL + \bM)^{-1}(\bL - \bM) \right)$ and converges monotonically in the norm $\|\cdot\|_{\bP}$ to the unique critical point. 
    Finally, notice that in the case $\bL-\bM$ is singular, the infimal convergence rate may not be attained, since $r$ can go towards infinity if $\rho$ goes towards $\rho_\mathrm{grad}$. However, if $\bL-\bM$ is non-singular, then $r$ disappears from the equation and the constructed gradient method converges at the rate $\rho_\mathrm{grad}$.
    %Finally, notice that for the case $\bM \prec_c\bL$ we have $\rho = \rho_\mathrm{grad}$.
\end{remark}

\begin{theorem}[Existence and uniqueness of critical points]
	\label{thm:ExMinCoCa}
	Let $\bM, \bL \in \bbR^{d\times d}$ be given symmetric matrices. Then the following three statements are equivalent:
	\begin{enumerate}
		\item The matrices $\bM$, $\bL$ are non-singular and satisfy $\bM \preceq_c \bL$.
		\item $S(\bM,\bL)$ is not empty and for all $f \in S(\bM,\bL)$ there exists at least one $z_f^* \in \bbR^d$ with $\nabla f(z_f^*) = 0$.
		\item $S(\bM,\bL)$ is not empty and for all $f \in S(\bM,\bL)$ there exists at most one $z_f^* \in \bbR^d$ with $\nabla f(z_f^*) = 0$.
		%\item $S(\bM,\bL)$ is not empty and for all $f \in S(\bM,\bL)$ exists at 
		%exactly one $z_f^* \in \bbR^d$ with $\nabla f(z_f^*) = 0$.
%		\item $S(\bM,\bL)$ is not empty and for all $f \in S(\bM,\bL)$ exists at most one $z_f^* \in \bbR^d$ with $\nabla f(z_f^*) = 0$.
	\end{enumerate}
\end{theorem}

\begin{remark}
	Theorem \ref{thm:ExMinCoCa} shows that if we aim to design algorithms that are convergent for the whole class $S(\bM,\bL)$, we must necessarily require $\bM \preceq_c \bL$, because otherwise there would be elements of $S(\bM,\bL)$ without critical points. Hence the introduced partial ordering plays a key role in our results.
	Note that it is no coincidence that in Theorem \ref{thm:ExMinCoCa} the existence of critical points for all functions in $S(\bM,\bL)$ and the uniqueness of critical points are two separate, equivalent statements. Similar to solutions of linear equation systems, here a solution for the equation $\nabla f(z) = 0$ exists for all $f \in S(\bM,\bL)$ if and only if the solution is unique for all $f \in S(\bM,\bL)$.
\end{remark}

\section{Main results}

In this section, a convex synthesis approach of optimizer parameters $(\bA,\bB,\bC)$ for the set of objective functions $S(\bM,\bL)$ and for a given convergence rate is provided. By Theorem \ref{thm:ConvTrInv}, the design for the class $S(\bM,\bL)$ reduces to designing algorithms for the class $S_0(0,\widetilde{\bL}) = \{ f \in S(0,\widetilde{\bL})| \nabla f(0) = 0 \}$ with $\widetilde{\bL} = \bL - \bM$. Hence we consider Problem \ref{problem2} instead of Problem \ref{problem1}.

\subsection{A Class of Lyapunov functions}

To design the algorithm parameters $(\bA,\bB,\bC)$ with a predescribed convergence rate, we propose the following class of (non-quadratic) Lyapunov function candidates
\begin{align*}
V_f(x) =& \begin{pmatrix}
x\\
\nabla f(\bC x)
\end{pmatrix}^T
\begin{pmatrix}
\bP_{11} & \bP_{12}\\
\bP_{21} & \bP_{22}
\end{pmatrix}
\begin{pmatrix}
x\\
\nabla f(\bC x)
\end{pmatrix}\\
&+ f(\bC x) - f(0) - \frac{1}{2}\nabla f(\bC x)^T \widetilde{\bL}^\dagger\nabla f(\bC x)
\end{align*}
with parameter $0 \prec \bP = \bP^T \in \bbR^{n+d\times n+d}$. (Recall, that in Problem \ref{problem2}, $\widetilde{\bL}$ was defined as $\bL - \bM$.) %Similar Lyapunov functions have been proposed by Yakubovich for the case $d = 1$ and continuous time systems in \cite{yakubovich1965method} and have been employed e.g. in \cite{josselson1974absolute} and \cite{suykens1998absolute}. The similarity is the first term in the Lyapunov function, which is quadratic term in the state and the nonlinearity.
%\ce{Check paper: Stability Criteria of Sector- and Slope-Restricted Lur’e Systems PooGyeon Park and follow up work. Julien M. Hendrickx} \dg{This paper uses similar Lyapunov functions.}
%\ce{Discuss similarities and differences in more detail in the paper. Write a few sentences about this and related Lyapunov functions.}
Similar Lyapunov functions have already been applied to Lur'e systems in continuous-time. Those Lyapunov functions share the first term, which is quadratic in the state $x$ and the static non-linearity $\nabla f(z)$. They have been proposed by Yakubovic for the case $d = 1$ in \cite{yakubovich1965method} and have been employed e.g. in \cite{josselson1974absolute}, \cite{park2002stability} and \cite{suykens1998absolute}. 

Our design approach, for a given convergence rate $\rho$, is based on finding simultaneously a Lyapunov function ($\bP \succ 0$) and algorithm parameters 
by semi-definite programming such that the Lyapunov conditions in the next theorem are satisfied.

\begin{theorem}[Lyapunov function and convergence rate for the algorithms]
	\label{thm:2.1LyapOpt}
	Let $(\bA,\bB,\bC)$ be parameters of Algorithm \eqref{eq:algo} for the set of objective functions $S$. If there exists a family of function $V_f: \bbR^n \to [0,\infty [$ satisfying quadratic bounds
	\begin{align}
	\alpha_f \| x - x_f^*\|^2 &\leq V_f(x) \leq \beta_f \| x - x_f^*\|^2 & \forall x \in \bbR^n, f\in S \label{eq:2.1LyaOpt}
	\end{align}
	for some fixed $\alpha_f,\beta_f \in \bbR_{>0}$ and the $\rho$-weighted increment bound
	\begin{align}
	V_f(x^+) - \rho^2 V_f(x) &\leq 0 & \forall x \in \bbR^d, \label{eq:2.1LyaOpt2}
	\end{align}
	where $x^+ = \bA x + \bB \nabla f(\bC x)$, then the optimizer defined by \eqref{eq:algo} is convergent with rate $\rho$.
\end{theorem}

The following two lemmas provide useful bounds for the considered class of Lyapunov functions and their increments and imply as by-product the positive definiteness of $V_f$.

\begin{lemma}[Quadratic bounds on $V_f$]
	\label{lem:4.4QuadBoundVf}
	Let $f \in S_0(0,\widetilde{\bL})$. Then the Lyapunov function candidates $V_f$ fulfill the quadratic bounds
	\begin{align*}
	\alpha_f \| x\|^2 \leq V_f(x) \leq \beta_f \|x\|^2
	\end{align*}
	with the constants $\alpha_f := \lambda_\mathrm{min}(\bP)$ and $\beta_f := \lambda_\mathrm{max}(\bP)(1+\|\widetilde{\bL}\|^2\|\bC\|^2) + \frac{\|\widetilde{\bL}\|\|\bC\|^2}{2}$.
\end{lemma}

\begin{lemma}[Upper bound on the Lyapunov increment of $V_f$]
	\label{lem:4.4EstDeriv}
	Assume $f \in S_0(0,\widetilde{\bL})$. Then, the weighted increment of $V_f$ from \eqref{eq:2.1LyaOpt2} is upper bounded as follows:
	\begin{align*}
	&V_f(x^+) - \rho^2 V_f(x) \leq\\
	&\resizebox{\linewidth}{!}{$
	\left(
	\begin{array}{c}
	x\\
	w\\
	\hline
	x^+\\
	w^+
	\end{array}\right)^T
	\left(
	\begin{array}{cc|cc}
	-\rho^2 \bP_{11} & -\rho^2 \bP_{12} & 0 & 0\\
	-\rho^2 \bP_{21} & -\rho^2 \bP_{22} & 0 & 0\\
	\hline
	0 & 0 & \bP_{11} & \bP_{12}\\
	0 & 0 & \bP_{21} & \bP_{22}
	\end{array}
	\right)
	\left(
	\begin{array}{c}
	x\\
	w\\
	\hline
	x^+\\
	w^+
	\end{array}\right)$}\\
	&\resizebox{\linewidth}{!}{$+
	\left(
	\begin{array}{c}
	x\\
	w\\
	\hline
	x^+\\
	w^+
	\end{array}\right)^T
	\left(
	\begin{array}{cc|cc}
	0 & 0 & 0 & -\frac{\lambda}{2}\bC^T\\
	0 & 0 & 0 & \frac{\lambda}{2}\widetilde{\bL}^\dagger\\
	\hline
	0 & 0 & 0 & \frac{1}{2}\bC^T\\
	-\frac{\lambda}{2}\bC & \frac{\lambda}{2}\widetilde{\bL}^\dagger & \frac{1}{2}\bC & -\widetilde{\bL}^\dagger
	\end{array}
	\right)
	\left(
	\begin{array}{c}
	x\\
	w\\
	\hline
	x^+\\
	w^+
	\end{array}\right),$}
	\end{align*}
	where $w = \nabla f(\bC x)$, $w^+ = \nabla f(\bC x^+)$ and $x^+ = \bA x + \bB w$, for arbitrary $\lambda \in [0,\rho^2]$.
\end{lemma}

\subsection{Convex synthesis of algorithms}

The following theorem reformulates the condition \eqref{eq:2.1LyaOpt2} in Theorem \ref{thm:2.1LyapOpt} using the established bound in Lemma \ref{lem:4.4EstDeriv}.

\begin{theorem}[Analysis Inequalities]
	\label{thm:4.5}
	Let $\bA \in \bbR^{n\times n}, \bB \in \bbR^{n\times d}$ and $\bC \in \bbR^{d\times n}$ be given. Set $\widetilde{\bA} = \bA + \bB\bM\bC$. Then the gradient-based algorithm \eqref{eq:algo} solves Problem \ref{problem1} and has convergence rate $\rho \in [0,1[$, if there exist $\bP = \bP^T \succ 0$, $\lambda \in [0,\rho^2]$ and $r \in \bbR$ such that the constraint \eqref{eq:Constraint}, i.e. $\bI_d = \bC (\widetilde{\bA} - \bI)^{-1} \bB \bM$, is satisfied and
	{\small
	\begin{align}
	&\resizebox{\linewidth}{!}{$
	\left(
		\begin{array}{ccc}
		\bI_n & 0 & 0\\
		0 & \bI_d & 0\\
		\hline
		\widetilde{\bA} & \bB & 0\\
		0 & 0 & \bI_d
		\end{array}\right)^T
	\left(
	\begin{array}{cc|cc}
	-\rho^2 \bP_{11} & -\rho^2 \bP_{12} & 0 & 0\\
	-\rho^2 \bP_{21} & -\rho^2 \bP_{22} & 0 & 0\\
	\hline
	0 & 0 & \bP_{11} & \bP_{12}\\
	0 & 0 & \bP_{21} & \bP_{22}
	\end{array}
	\right)
	(\star) $}\nonumber\\
	&\resizebox{\linewidth}{!}{$
	+
	\left(
		\begin{array}{ccc}
		\bI_n & 0 & 0\\
		0 & \bI_d & 0\\
		\hline
		\widetilde{\bA} & \bB & 0\\
		0 & 0 & \bI_d
		\end{array}\right)^T
	\left(
	\begin{array}{cc|cc}
	0 & 0 & 0 & -\frac{\lambda}{2}\bC^T\\
	0 & -r\bPi & 0 & \frac{\lambda}{2}\widetilde{\bL}^\dagger\\
	\hline
	0 & 0 & 0 & \frac{1}{2}\bC^T\\
	-\frac{\lambda}{2}\bC & \frac{\lambda}{2}\widetilde{\bL}^\dagger & \frac{1}{2}\bC & -\widetilde{\bL}^\dagger-r\bPi
	\end{array}\right)
	(\star)$}\nonumber \\
	& \hspace{73mm} \prec 0 \label{eq:4.5_1}
	\end{align}}
	is satisfied, where $\widetilde{\bL} = \bL-\bM$ and $\bPi = \bPi_{\ker (\bL-\bM)}$.
\end{theorem}

Theorem \ref{thm:4.5} provides sufficient conditions for a given algorithm to achieve a convergence rate $\rho$. Notice that the conditions in Theorem \ref{thm:4.5}
are affine in the positive definite decision variable $\bP$ and hence semi-definite programming can be used to verify these conditions.
For the synthesis of algorithms, i.e. if in addition to $\bP$ also $\bA,\bB,\bC$ are decision
variables, the decision variables enter in a non-affine (non-convex) fashion and thus,
an efficient synthesis of algorithms with semi-define programming is not possible.
Hence, it is of key importance to find \emph{equivalent} conditions in terms of matrix inequalities and equations in which the decision variables enter in an affine fashion.
The following theorem shows that this is indeed possible.

\begin{theorem}[Synthesis Inequalities]
	\label{thm:4.6SynIneq}
	Let $n\geq 3d$. Then there exist matrices $\bA\in \bbR^{n\times n}, \bB \in \bbR^{n\times d}, \bC \in \bbR^{d\times n}$, which render the conditions
	\eqref{eq:4.5_1} and \eqref{eq:Constraint} in Theorem \ref{thm:4.5}
	for a given convergence rate $\rho$   feasible, if and only if there exist $\hat{\bA} \in \bbR^{n\times n}, \hat{\bB} \in \bbR^{n\times d}, \bC \in \bbR^{d\times n}, \bP = \bP^T \in \bbR^{n+d\times n+d}$, $r \in \bbR$ and $\lambda \in [0,\rho^2]$ such that the matrix inequality
	\begin{align}
	&\resizebox{\linewidth}{!}{$
	\left(
	\begin{array}{cc|c|cc}
	-\rho^2 \bP_{11} & -\rho^2 \bP_{12} & * & * & * \\
	-\rho^2 \bP_{21} & -\rho^2 \bP_{22}-r\bPi & * & * & * \\
	\hline
	\frac{1}{2}\bJ_2\hat{\bA}  -\frac{\lambda}{2}\bC & \frac{1}{2} \bJ_2\hat{\bB}  + \frac{\lambda}{2}\widetilde{\bL}^\dagger & -\widetilde{\bL}^\dagger-r\bPi & * & *\\
	\hline
	\hat{\bA} & \hat{\bB} & \bP_{12} & -\bP_{11} & -\bP_{12}\\
	\bJ_3\hat{\bA} & \bJ_3\hat{\bB} & \bP_{22} & -\bP_{21} & -\bP_{22}
	\end{array}
	\right)$}\nonumber\\
	& \hspace{71mm} \prec 0,	\label{eq:SynIneq}
	\end{align}
	with $\widetilde{\bL} = \bL-\bM$ and $\bPi = \bPi_{\ker (\bL-\bM)}$
	is satisfied and the constraints
	\begin{align}
    	\begin{split}
    	\hat{\bB} &= (\hat{\bA} - \bP_{11})\bJ_1^T\bM^{-1}, \hspace{10mm} \bC \bJ_1^T = \bI_d,\\ \bC &= \bJ_2 \bP_{11}, \hspace{24mm} \bP_{21} = \bJ_3 \bP_{11},
    	\end{split} \label{eq:ConstraintSyn}
	\end{align}
	are satisfied, where $\bJ_1,\bJ_2,\bJ_3 \in \bbR^{d\times n}$ are
	{\small
	\begin{align*}
	\resizebox{\linewidth}{!}{$
	\bJ_1 = \begin{pmatrix}
	\bI_d & 0
	\end{pmatrix},
	\hspace{2mm}
	\bJ_2 = \begin{pmatrix}
	0_d & \bI_d & 0
	\end{pmatrix},
	\hspace{2mm}
	\bJ_3 = \begin{pmatrix}
	0_d & 0_d & \bI_d & 0
	\end{pmatrix}.$}
	\end{align*}}
	The algorithm parameters $\bA,\bB,\bC$ can be obtained by
	\begin{align*}
	\bB = \bP_{11}^{-1}\hat{\bB}, \hspace{10mm}
	\bA = \bP_{11}^{-1}\hat{\bA}-\bB\bM\bC.
	\end{align*}
\end{theorem}

Finally, we want to show that the above matrix inequalities are always feasible,
by analyzing or designing gradient algorithms. %\ce{Skip theorem and proof.}

\begin{theorem}[Existence of Solutions]
	\label{thm:Existence of solutions}
	The following four statements are equivalent:
	\begin{enumerate}[i)]
		\item The matrices $\bM,\bL$ are non-singular and satisfy $\bM \preceq_c \bL$.
		\item The gradient method defined by $(\bA,\bB,\bC)$ with
		$
			\bA = \bI_d, \hspace{2mm} \bB = -2(\bL+\bM)^{-1}, \hspace{2mm} \bC = \bI_d
		$
		fulfills the conditions \eqref{eq:4.5_1} and \eqref{eq:Constraint} of Theorem \ref{thm:4.5} for any $\rho \in ]\rho_\mathrm{grad},1[$, where
			$\rho_\mathrm{grad} = \rho\left( (\bL + \bM)^{-1} (\bL - \bM) \right).$
		\item For all $n\geq d$, there exists an algorithm $(\bA,\bB,\bC)$, $f \in S(\bM,\bL)$ and $\rho \in [0,1[$ such that the conditions
		\eqref{eq:4.5_1} and \eqref{eq:Constraint}
		in Theorem \ref{thm:4.5}  are satisfied.
		\item For all $n \geq 3d$ there exists a solution to \eqref{eq:SynIneq} and \eqref{eq:ConstraintSyn} in Theorem \ref{thm:4.6SynIneq} for some $\rho \in [0,1[$.
	\end{enumerate}
\end{theorem}

If one optimizes simultaneously over $\bA,\bB,\bC$ and $\rho$ in Theorem \ref{thm:4.5} or Theorem \ref{thm:4.6SynIneq}, then $\rho_\mathrm{grad}$ is usually not the optimal rate. Often, there exist faster algorithms. However, it is the optimal rate for the gradient method from Proposition~\ref{prop:4.1GradContraction}.

\subsection{Comparison to IQC based approaches}

As already mentioned in the introduction, there exist quite some publications on the application of methods from robust control theory to gradient-based optimization. 
Some of these approaches use a technique called IQCs (integral quadratic constraints). The goal of this subsection is to explain the relation between IQC based approaches 
(such as in \cite{lessard2016analysis}, \cite{michalowsky2019robust})
and the Lyapunov based approach as in this paper. For this purpose, we will restrict ourselves to the special case of $S(m,l)$ with scalar $m,l \in \bbR_{\geq 0}$, which is usually considered in the literature.
%A tool for the analysis of exponential stability are so-called $\rho$-IQCs,
%see for example \cite{lessard2016analysis} or \cite{boczar2017exponential}.
%\begin{definition}[$\rho$-IQC]
	%\label{def:rhoIQC}
	
The main steps of IQC based approaches are summarized in the subsequent paragraphs:
\begin{enumerate}[(a)]
    \item Let $y \in \ell_{2,\rho}$ and $u\in \ell_{2,\rho}$ be signals with associated $z$-transforms $\hat{y}(z)$ and $\hat{u}(z)$.
	Then these signals are said to satisfy the $\rho$-IQC defined by a Hermitian complex-valued function $\bPi$ if
% 	\begin{align}
% 	\int_{\bbC_{|z| = \rho}} \begin{pmatrix}
% 	\hat{y}(z)\\
% 	\hat{u}(z)
% 	\end{pmatrix}^*
% 	\bPi(z) \begin{pmatrix}
% 	\hat{y}(z)\\
% 	\hat{u}(z)
% 	\end{pmatrix}dz \geq 0. \label{eq:3.1IQC_1}
% 	\end{align}
	\begin{align}
	\int_{0}^{2\pi} \begin{pmatrix}
	\hat{y}(\rho e^{i\omega})\\
	\hat{u}(\rho e^{i\omega})
	\end{pmatrix}^*
	\bPi(\rho e^{i\omega}) \begin{pmatrix}
	\hat{y}(\rho e^{i\omega})\\
	\hat{u}(\rho e^{i\omega})
	\end{pmatrix}\diff \omega \geq 0. \label{eq:3.1IQC_1}
	\end{align}
	A bounded causal operator $\Delta$ satisfies the $\rho$-IQC defined by $\bPi$ if \eqref{eq:3.1IQC_1} holds for all $y \in \ell_2$ and $u = \Delta (y)$. $\mathrm{IQC}(\bPi,\rho)$ denotes the set of all $\Delta$ that satisfy the $\rho$-IQC defined by $\bPi$. 
	\item Next view the gradient-based algorithm from Problem \ref{problem2} as an interconnection of the linear system defined by the transfer function $\bG(z) = \bC (z\bI_n - \widetilde{\bA})^{-1}\bB$ and the static nonlinearity defined by $(y_k)_{k\in\bbN_0} \mapsto (\nabla f(y_k))_{k\in\bbN_0} =:\Delta(y)$. It is well-known (see for example \cite{lessard2016analysis}) that $\nabla f$ satisfies the IQC, i.e.
	an operator (system) $\Delta: \ell_{2e}^d \to \ell_{2e}^d$ which is static and slope restricted in the sector $[m,l]$ satisfies the IQC defined by the multiplier
	\begin{align}
	    \bPi(z) = \bPsi^*(z)\bR \bPsi(z), 
	    \label{offbyone}
	\end{align}
	where the factorization is given by
	\begin{align*}
        \resizebox{\linewidth}{!}{$
        \bPsi(z) = \begin{pmatrix}
    	(l-m)(1-\lambda z^{-1}) & z^{-1}\lambda \bI_d\\
    	0 & \bI_d
    	\end{pmatrix},
        \bR = \begin{pmatrix}
        0 & \bI_d\\
        \bI_d & -2 \bI_d
        \end{pmatrix}.$}
    \end{align*}
% 	\begin{align}
% 	    \resizebox{\linewidth}{!}{$
% 	    \begin{pmatrix}
% 	        -ml\left(2 - \lambda z^{-1} - \lambda\bar{z}^{-1}\right) & m ( 1- \lambda z^{-1}) + l(1 - \lambda\bar{z}^{-1})\\
% 	        m ( 1- \lambda \bar{z}^{-1}) + l(1 - \lambda z^{-1}) & -\left(2 - \lambda z^{-1} - \lambda\bar{z}^{-1}\right) 
% 	    \end{pmatrix}.$}
% 	\end{align}
	\item Finally, the following result from IQC theory (see e.g. \cite{boczar2017exponential,michalowsky2019robust}), which is based on the exponential weighting operators $\rho_{\pm}: \ell_{2e}(\bbN_0) \to \ell_{2e}(\bbN_0), (u_k) \mapsto (\rho^{\pm k} u_k)$, is invoked 
    to verify convergence of the algorithm with rate~$\rho$.
    \begin{theorem}[Exponential stability with IQCs]
    	\label{thm:3.1ExpStabIQC}
    	Fix $\rho \in ]0,1[$. Let $G$ be a stable, causal linear dynamical system with transfer function $\bG(\cdot)$ and all poles of $\bG$ contained in $\bbC_{|z| < \rho}$. Let further $\Delta$ be a stable, causal dynamical system such that $\Delta' = \rho_- \circ \Delta \circ \rho_+$ is a bounded operator. Suppose that:
    	\begin{enumerate}[i)]
    		\item for all $\tau \in [0,1]$, the interconnection of $G$ and $\tau \Delta$ is well posed,
    		\item for all $\tau \in [0,1]$, we have $\tau \Delta \in \mathrm{IQC} (\bPi,\rho)$,
    		\item the following frequency domain inequality (FDI) holds:
    		\begin{align}
        		\begin{pmatrix}
        		\bG(z)\\
        		\bI
        		\end{pmatrix}^*
        		\bPi(z)
        		\begin{pmatrix}
        		\bG(z)\\
        		\bI
        		\end{pmatrix}&\prec 0, & \forall z \in \bbC_{|z| = \rho}. \label{eq:FDI}
    		\end{align}
    	\end{enumerate}
    	Then, the feedback interconnection of $G$ and $\Delta$ is exponentially stable with rate $\rho$.
    \end{theorem}
\end{enumerate}

The connection to the proposed Lyapunov based approach can now be established by applying Theorem~\ref{thm:3.1ExpStabIQC} to Problem~\ref{problem2}:
\begin{itemize}
    \item The interconnection of $\bG$ and $\nabla f$ is always well posed, because $\bG$ is strictly proper and $\nabla f$ has relative degree zero.
    \item Condition $ii)$ is always satisfied for $\nabla f$ with $f \in S_0(0,l\bI_d - m \bI_d)$. (Note, that $ii)$ is not satisfied for $f \in S(m\bI_d,l\bI_d)$.)
    \item Finally, boundedness of $\Delta'$ is a consequence of the Lipschitz continuity of $\nabla f$.
\end{itemize}

The following Lemma  states the relation between the frequency domain inequality \eqref{eq:FDI} in iii) and the matrix inequality \eqref{eq:4.5_1} from Theorem \ref{thm:4.5}.

\begin{lemma}[Relation between IQC and Lyapunov-based approach]
    \label{lem:FDI-LMI-connection}
    The FDI \eqref{eq:FDI} and $\sigma (\widetilde{\bA}) \subsEq \bbC_{|z| < \rho}$ hold if and only if \eqref{eq:4.5_1} is feasible for some $\lambda \in [0,\rho^2]$.
\end{lemma}

All together, it is possible to prove Theorem \ref{thm:4.5} using Theorem \ref{thm:3.1ExpStabIQC} and Lemma \ref{lem:FDI-LMI-connection}.

\section{Examples and numerical results}
\label{sec:5}

%Here, a few numerical results on the optimizers from this paper are shown. The first result is on the performance of methods synthesized by solving the matrix inequality from Theorem \ref{thm:4.6SynIneq}.

\subsection{Convergence rates}

\begin{figure}
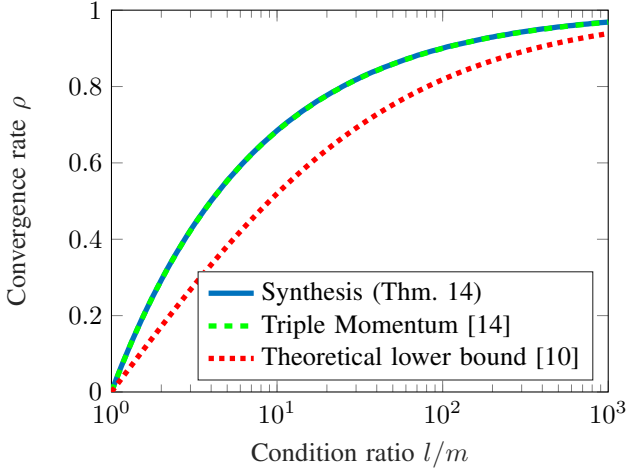

	\centering
	\include{PerformancePaper}
	\caption{The convergence rate guarantees achieved by designing algorithms for $S(m,l)$ using Theorem \ref{thm:4.6SynIneq} are plotted over the condition number $l/m$ and compared to the rate bound of Triple Momentum $\rho = 1 - \frac{\sqrt{m}}{\sqrt{l}}$ and the theoretical lower bound $\rho = \frac{\sqrt{l} - \sqrt{m}}{\sqrt{l}+\sqrt{m}}$ from \cite{nesterov1998introductory}.}
	\label{fig:6.1ConvRates}
\end{figure}

To demonstrate the performance of our synthesis, we apply it to the class $S(m,l)$ of strongly convex functions, which is often considered in the literature (for example in \cite{lessard2016analysis}, \cite{michalowsky2019robust} and \cite{nesterov1998introductory}). The algorithm parameters $(\bA,\bB,\bC)$ are computed by solving \eqref{eq:SynIneq} and \eqref{eq:ConstraintSyn} in Theorem \ref{thm:4.6SynIneq} for $\lambda = \rho^2$, where $\rho$ is optimized using a bisection search. Here, setting $\lambda$ equal to $\rho^2$ is motivated by the proof of Lemma~\ref{lem:4.4EstDeriv}, where $\lambda = \rho^2$ gives the sharpest estimate on the increment of the Lyapunov function. The obtained convergence rates are shown in Figure~\ref{fig:6.1ConvRates}, where they are compared to the convergence rates of the Triple Momentum method from \cite{van2017fastest} and the theoretical lower bound on the convergence rates obtained by Nesterov. As can be observed, our synthesized algorithm has the same convergence rates as the Triple Momentum method. A result, that is also obtained in \cite{lessard2019direct} using an IQC based approach.

Strictly speaking, the synthesis with Theorem~\ref{thm:4.6SynIneq} is  not an LMI synthesis if we consider $\lambda$
as a decision variable. This parameter could possibly be optimized using a line search algorithm, which we did in the first place. However, in our empirical experiments, we found that in the case $\bM \succ 0$ the value $\lambda = \rho^2$ was always the optimal one.

\subsection{Structured objective functions}

\begin{figure}
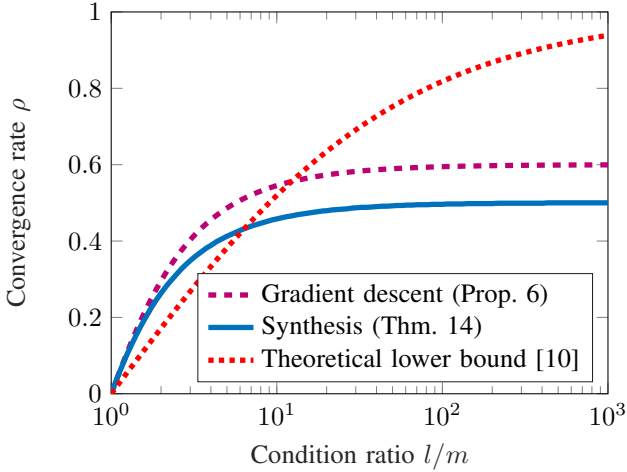

	\centering
	\include{benefitsStructuredMethods}
	\caption{Convergence rates achieved by the gradient descent algorithm in Proposition
		\ref{prop:4.1GradContraction}
	and by synthesis with Theorem \ref{thm:4.6SynIneq} for $S(\bM,\bL)$.
	Note that the theoretical lower bound holds for the class $S(m,l)$ and not for the subset $S(\bM,\bL) \subset S(m,l)$, because the subset contains fewer objective functions. \label{fig:structuredGD}}
\end{figure}

The following (academic) example shall demonstrate the possible benefits of including additional properties of the objective function into algorithm design compared to the design for $S(m,l)$. Consider the class of functions $S(\bM,\bL)$ with
\begin{align*}
	\bM &= \begin{pmatrix}
	l-m+\frac{m^2}{l} & 0\\
	0 & m
	\end{pmatrix},~~
	\bL = \bS^T
	\begin{pmatrix}
	l & 0\\
	0 & 2m-\frac{m^2}{l}
	\end{pmatrix}
	\bS,\\
	\bS &= 
	\begin{pmatrix}
	\sqrt{1-\left(\frac{m}{l}\right)^2} & -\frac{m}{l}\\
	\frac{m}{l} & \sqrt{1-\left(\frac{m}{l}\right)^2}
	\end{pmatrix}.
\end{align*}
These matrices fulfill $m\bI \preceq \bM \preceq_c \bL \preceq l \bI$. Moreover, the largest eigenvalue of $\bL$ is $l$ and the smallest eigenvalue of $\bM$ is $m$. Hence, the best \quot{standard method} for the class $S(\bM,\bL)$ is a method for $S(m,l)$ and has a convergence rate that is not faster than $\frac{\sqrt{l} - \sqrt{m}}{\sqrt{l} + \sqrt{m}}$. The method designed using Theorem \ref{thm:4.6SynIneq} on the other hand has at least the convergence rate $\rho((\bM + \bL)^{-1}(\bL-\bM))$. Figure \ref{fig:structuredGD} illustrates these convergence rates together with the rate of a synthesized algorithm. One can recognize that in this example the structured method is superior to any unstructured method.
% These matrices fulfil $m\bI \preceq \bM \preceq_c \bL \preceq l \bI$. Moreover, the largest eigenvalue of $\bL$ is $l$ and the smallest eigenvalue of $\bM$ is $m$. Hence, any method designed for $S(m,l)$, which does not take the additional information given by $\bM$ and $\bL$ into account has a convergence rate no better than $\frac{\sqrt{l} - \sqrt{m}}{\sqrt{l} + \sqrt{m}}$. The method designed for $S(\bM,\bL)$ using Theorem \ref{thm:4.6SynIneq} on the other hand has at least the convergence rate $\rho((\bM + \bL)^{-1}(\bL-\bM))$. Figure \ref{fig:structuredGD} illustrates these convergence rates together with the rate of a synthesized algorithm. One can recognize that in this (academic) case the structured methods is superior to any unstructured method.

\subsection{Application to constrained optimization}

The class $S(\bM,\bL)$ can contain non-convex functions. If both $\bM$ and $\bL$ are indefinite but the condition $\bM \preceq_c \bL$ is fulfilled, then $S(\bM,\bL)$ is a class of functions with unique critical (saddle) points. 
One particular saddle point problem can be obtained in the context of convex constrained optimization. If one aims to solve the (linearly) constrained optimization problem
\begin{align}
\begin{split}
	\mathrm{minimize}~ & f(x), \label{eq:ConstraintOpt} \\
	                                 \mathrm{subject~to~} & x\in \bbR^d,~  \bA_\mathrm{eq} x = b_\mathrm{eq},
\end{split}	                                 
\end{align}
where $f\in S(\bM,\bL)$, $\bA_\mathrm{eq} \in \mathbb{R}^{d_2 \times d}$ and $0 \prec \bM \prec \bL$ holds (such that $f$ is strictly convex), then a solution can be found by solving the saddle point  problem
\begin{align*}
	\sup_{\lambda \in \bbR^{d_2}} \inf_{x\in \bbR^d} f(x) + \lambda^T(\bA_\mathrm{eq} x - b_\mathrm{eq}).
\end{align*}
Here, the Lagrangian function 
$L(x,\lambda)=f(x) + \lambda^T(\bA_\mathrm{eq} x - b_\mathrm{eq})$
is an element of $S(\bM_L,\bL_L)$, where
\begin{align}
	\bM_L = \begin{pmatrix}
	\bM & \bA_\mathrm{eq}^T\\
	\bA_\mathrm{eq} & 0
	\end{pmatrix}, \hspace{3mm} \bL_L = \begin{pmatrix}
	\bL & \bA_\mathrm{eq}^T\\
	\bA_\mathrm{eq} & 0
	\end{pmatrix}.
	\label{tmtl}
\end{align}
If $\bM_L \preceq_c \bL_L$ is satisfied, then our design procedure can be applied to design a gradient based algorithm for $L$, which solves the constrained optimization problem. The following lemma shows under rather mild conditions that this is possible. 
\begin{lemma}
    \label{lem:5.1}
	Let $\bA_\mathrm{eq} \in \bbR^{d_2 \times d}$ and symmetric matrices $\bM,\bL \in \bbR^{d\times d}$ be given.
	 Consider $\bM_L,\bL_L$ defined in \eqref{tmtl} and
	assume that $\bM \preceq_c \bL$ holds with $\bM,\bL$ being non-singular and that $\bA_\mathrm{eq}$ has full row rank. Then
$
	\bM_L  \preceq_c \bL_L
$
	holds, and $\bM_L$ and $\bL_L$ are non-singular.
\end{lemma}

\begin{figure}
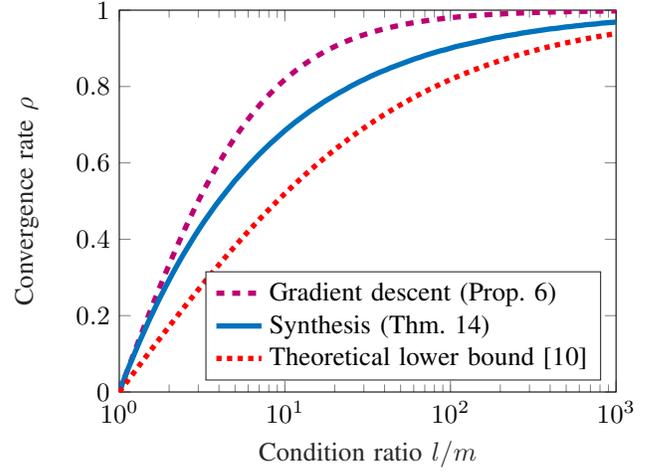

	\centering
	\include{ConstrainedOpt}
	\caption{Convergence rates achieved by the gradient descent
	algorithm  in Proposition \ref{prop:4.1GradContraction} and
	by designing algorithms with Theorem \ref{thm:4.6SynIneq} for a constrained optimization problem. \label{fig:ConstrainedRates}}
\end{figure}

As an academic example, consider the constrained optimization problem \eqref{eq:ConstraintOpt} with $f \in S(m\bI_2,l\bI_2)$ and $\bA_\mathrm{eq} = (1 \hspace{2mm} 1)$. 
As described above, matrices $\bM_L \preceq_c \bL_L$ can be constructed such that the Lagrangian $L$ of this problem is in $S(\bM_L,\bL_L)$. This enables algorithms of the form
$
x_{k+1}
=
\bA
x_k + 
\bB
\nabla L (\bC x_k)$ to be designed.
The algorithm parameters $\bA,\bB,\bC$ can be designed by solving the matrix inequality from Theorem \ref{thm:4.6SynIneq}. 	
The results are presented in Figure \ref{fig:ConstrainedRates}.
For the sake of comparison, we added the rates of the descent algorithm from Proposition \ref{prop:4.1GradContraction}.
Interestingly, the convergence rates are exactly equal to the convergence rates for the unconstrained optimization problems. 
In general, we observed in our experiments that the convergence rates for linearly constrained optimization problems were often faster than those for unconstrained problems, but never slower.
Notice that we have the condition $n \ge 3d$ in Theorem \ref{thm:4.6SynIneq}, hence
the algorithm with one equality constraint has at least dimension 9. However,
it is often possible to reduce the dimension of the algorithm as outlined below.
For example, we consider the algorithm parameters $\bA,\bB,\bC$ for $m = 1,l=15$ designed using Theorem \ref{thm:4.6SynIneq}. The original matrices had dimension $n = 9$. 
We observed that the last three modes usually do not contribute much to the dynamics of the algorithm. Hence, it is possible to eliminate them using balanced truncation. 
In our example, we obtained the reduced parameters: 
\begin{align*}
    \bA &= \begin{pmatrix}
    1 & 0 & 0 &   0.0135  & -0.0258  & -0.0017\\
    0 & 1 & 0 &   0.0135  &  0.0258  & -0.0017\\
    0 & 0 & 1 &  -0.6076  & -0.0036  & -0.0363\\
    0 & 0 & 0 &  -0.3097  & -0.0042  & -0.0474\\
    0 & 0 & 0 &  -0.0039  &  0.3909  & -0.0002\\
    0 & 0 & 0 &   1.1631  &  0.0070  &  0.5255
    \end{pmatrix},\\
    \bB &=\begin{pmatrix}
    -0.0846  &  0.0707  & -0.1978\\
     0.0707  & -0.0846  & -0.1978\\
    -0.2758  & -0.2758  & -3.2399\\
     0.0860  &  0.0940  & -4.7039\\
     0.6738  & -0.6727  & -0.0264\\
     0.0896  &  0.0900  &  6.3240\\
    \end{pmatrix},\\
    \bC &=\begin{pmatrix}
    1 & 0 & 0 & 0 & 0 & 0\\
    0 & 1 & 0 & 0 & 0 & 0\\
    0 & 0 & 1 & 0 & 0 & 0
    \end{pmatrix}.
\end{align*}
We used Theorem \ref{thm:4.5} to check that the reduced algorithm still converges for $S(m = 1,l=15)$.
The reduced algorithm achieves a convergence rate of at least $0.7422$. which is faster than the rate of gradient descent, which is $0.8750$ and exactly as fast as the unreduced algorithm.
In this example, we have chosen a specific representation in which the first $d$ columns of $\bA$ are the first $d$ unit vectors in $\bbR^n$ and $\bC$ takes the form of an identity matrix concatenated with a zero block. The existence of such a representation is guaranteed by \eqref{eq:Constraint}. From this specific representation, it can be extracted that $\bA$ will always have $d$ eigenvalues at one. The modes with one eigenvalues play the role of a memory for the current best guess of the optimization algorithm and are therefore necessary.

\section{Conclusion}

We presented a \emph{convex} synthesis procedure to design gradient-based algorithms based on a general class of Lur'e Lyapunov functions and linear matrix inequalities. The class of objective functions, which was considered, generalizes the class of strongly convex functions and offers the possibility to incorporate additional information into the algorithm design. It should be emphasized that this class of functions also includes non-convex functions - in particular functions with saddle points. The usefulness of our novel function class was demonstrated, firstly, by showing that additional information about the objective function can boost the convergence rate of algorithms considerably and, secondly, by showing that it can be used to design algorithms for solving optimization problems with linear equality constraints.

Open future research questions are for example the design of distributed algorithms or the design of optimization
algorithms for problems with inequality constraints.
%We hope that the proposed design procedure for accelerated gradient methods for the generalized class of objective functions may be applied in other situations beyond saddle point problems.
%Moreover,, \ce{LinAlg ... solving linear systems of equations as suggested by DG}

\section{Acknowledgements}
This work was funded by Deutsche Forschungsgemeinschaft (DFG, German Research Foundation) under Germany’s Excellence Strategy - EXC 2075 -390740016

\bibliographystyle{plain}
\bibliography{sources}

\appendix

\section{Proofs}

\subsection{Projections and pseudo inverses}

The pseudo inverse $\bL^\dagger$ and projection matrix $\bPi_{\ker \bL}$/$\bPi_{\im \bL}$ onto the kernel/image of a symmetric matrix $\bL$ are used at several places in the proofs of this paper. 
Hence, some important formulas are summarized below.
%The key to understand these rules is the singular value decomposition of a matrix $\bA \in \bbR^{d\times d}$:
%\begin{align*}
%    \bA = \underbrace{
%    \left(
%	\begin{array}{c}
%	\bU_1\\
%	\hline
%	\bU_2
%	\end{array}
%	\right)^T
 %   }_{\bU^T}
    %\underbrace{
    %\left(
	%\begin{array}{c|c}
	%\begin{pmatrix}
	%\sigma_1 &  & \\
%	 & \ddots & \\
	 %&  & \sigma_s
	%\end{pmatrix}
%	 & 0\\
	%\hline
	%0 & 0
	%\end{array}
	%\right)
    %}_{\bSigma}
    %\underbrace{
    %\left(
	%\begin{array}{c}
	%\bV_1\\
	%\hline
	%\bV_2
	%\end{array}
	%\right)
    %}_{\bV}
%\end{align*}
%with a diagonal matrix $\bSigma$ and orthogonal matrices $\bU,\bV$. This decomposition %is very insightful, as it allows us to express the pseudo inverse, the kernel projector and the image projector in a simple and intuitive way:
Let $\bA =  \bU^T \bSigma \bV$ be the singular value decomposition of a matrix $\bA$,
then
\begin{align*}
    \bA^\dagger &= \resizebox{0.83\linewidth}{!}{$\underbrace{
    \left(
	\begin{array}{c}
	\bV_1\\
	\hline
	\bV_2
	\end{array}
	\right)^T
    }_{\bV^T}
    \underbrace{
    \left(
	\begin{array}{c|c}
	\begin{pmatrix}
	\sigma_1^{-1} &  & \\
	 & \ddots & \\
	 &  & \sigma_s^{-1}
	\end{pmatrix}
	 & 0\\
	\hline
	0 & 0
	\end{array}
	\right)
    }_{\bSigma^\dagger}
    \underbrace{
    \left(
	\begin{array}{c}
	\bU_1\\
	\hline
	\bU_2
	\end{array}
	\right)
    }_{\bU}$}
\end{align*}
and    $\bPi_{\ker \bA} = \bV_2^T\bV_2$,     $\bPi_{\im \bA} = \bU_1\bU_1^T$.
%For the projectors, this is clear, for the Moore Penrose pseudo inverse, a proof can be %found in \cite{knabner2018lineare} (Theorem 4.129).\\
We will particularly be interested in the following four identities for the projectors and pseudo inverses of a symmetric positive semidefinite matrix $\bL$, $r \not = 0$:
%(in this case $\bU = \bV$):
\begin{align}
    \bPi_{\im \bL} &= \bL \bL^\dagger = \bL^\dagger \bL, \label{eq:PseudoInvProd}\\
    \bI_d &= \bPi_{\im \bL} + \bPi_{\ker \bL}, \label{eq:ProjektorSum}\\
    (\bL + r\bPi_{\ker \bL})^{-1} &= \bL^\dagger + \frac{1}{r}\bPi_{\ker \bL}, \label{eq:InverseWithPseudoInverse}\\
    (\bL^\dagger + r\bPi_{\ker \bL})^{-1} &= \bL + \frac{1}{r}\bPi_{\ker \bL}. \label{eq:InverseWithPseudoInverse2}
\end{align}
These identities follow from the singular value decomposition as shown above.
%using the above representations of the pseudo inverse is %standard linear algebra.

\subsection{Proof of Theorem \ref{thm:ConvTrInv}}

Step 1. First assume that $(\bA,\bB,\bC)$ solves Problem \ref{problem1}. We prove that $(\widetilde{\bA},\bB,\bC)$ solves Problem \ref{problem2}.
Let $g \in S_0(0,\widetilde{\bL}) = S_0(0,\bL-\bM)$ be an arbitrary function. Then
    \begin{align*}
 f(z) = g(z) + \frac{1}{2} z^T \bM z
    \end{align*}
    is an element of $S(\bM,\bL)$ with $\nabla f(0) = 0$. Now consider the iterates of algorithm \eqref{eq:algo} with the parameters $(\widetilde{\bA},\bB,\bC)$ for the objective functions $g$:
    \begin{align*}
        x_{k+1} &= \widetilde{\bA} x_k + \bB \nabla g(\bC x_k)\\
        &= \bA x_k + \bB (\nabla g(\bC x_k) + \bM\bC x_k)\\ % \underbrace{(\nabla g(\bC x_k) + \bM\bC x_k)}_{\nabla f(\bC x_k)}\\
        &= \bA x_k + \bB \nabla f(\bC x_k).
    \end{align*}
    Since those are the iterates of the algorithm defined by $(\bA,\bB,\bC)$ for $f \in S(\bM,\bL)$, we know that $x_k$ converges to $x_f^*$ at rate $\rho$ for any $x_0 \in \bbR^d$. Notice that $x_f^*$ must be zero in this case because zero is a fixed-point of the considered iteration (since $\bA 0 + \bB \nabla f(\bC 0) = 0$ by $\nabla f(0) = 0$) and hence, if $x_f^*$ were not zero, then the iterates for $x_0 = 0$ would not converge to $x_f^*$.
    It remains to show satisfaction of the constraint \eqref{eq:Constraint}.
    For this purpose define $f \in S(\bM,\bL)$ as
    \begin{align*}
        f(z) = \frac{1}{2} (z - z_f^*)^T\bM (z - z_f^*)
    \end{align*}
    for some $z_f^*\in \bbR^d$ and check that it satisfies $\nabla f(z_f^*) = 0$. By assumption, Problem \ref{problem1} is solved, meaning that the iterates
    \begin{align*}
        x_{k+1} &= \bA x_k + \bB \bM (\bC x_k - z_f^*)\\
        &= \widetilde{\bA}x_k - \bB \bM z_f^*
    \end{align*}
    of algorithm \eqref{eq:algo} converge to $x_f^*$ for any $x_0$. This implies, that $x_f^*$ is a solution of the fixed point equation
    \begin{align*}
        x_f^* = \widetilde{\bA} x_f^* - \bB \bM z_f^*.
    \end{align*}
    The convergence for arbitrary initial value implies that $\widetilde{\bA}$ is Schur and, hence, $\widetilde{\bA} - \bI_n$ must be non-singular. Then, the fixed point equation can be solved for $x_f^*$:
    \begin{align*}
        x_f^* = (\widetilde{\bA} - \bI_n)^{-1}\bB\bM z_f^*.
    \end{align*}
    By assumption, we have in addition
    \begin{align*}
        z_f^* = \bC x_f^* = \bC (\widetilde{\bA} - \bI_n)^{-1}\bB\bM z_f^*.
    \end{align*}
    Since $z_f^*$ is arbitrary, the constraint $\bC(\widetilde{\bA} - \bI_n)^{-1}\bB\bM = \bI_d$ must hold.
    
Step 2. Now assume, that $(\widetilde{\bA},\bB,\bC)$ is a solution of Problem \ref{problem2}. We prove that $(\bA,\bB,\bC)$ solves Problem \ref{problem1}.
For that, we first consider all functions $f \in S(\bM,\bL)$ for which there exists a critical point $z_f^*$.\\
Let $f \in S(\bM,\bL)$ be given such that there exists $z_f^*$ with $\nabla f(z_f^*) = 0$. Then $g$ defined by $g(z) = f(z + z_f^*) - \frac{1}{2}z^T\bM z$ is an element of $S_0(0,\bL - \bM) = S_0(0,\widetilde{\bL})$. Hence, the iterative scheme
\begin{align*}
    \tilde{x}_{k+1} = \widetilde{\bA} \tilde{x}_k + \bB \nabla g(\bC \tilde{x}_k)
\end{align*}
converges to zero at rate $\rho$ for any $x_0 \in \bbR^n$. Now add $x_f^* := (\widetilde{\bA} - \bI_n)^{-1}\bB\bM z_f^*$ on both sides of the above equation and consider the new sequence $x_k := \tilde{x}_k + x_f^*$:
\begin{align*}
    x_{k+1} &= \tilde{x}_{k+1} + x_f^*\\
    &= \widetilde{\bA}(\tilde{x}_k + x_f^*) + \bB \nabla g(\bC \tilde{x}_k) + x_f^* - \widetilde{\bA} x_f^*\\
    &= \bA x_k + \bB\bM\bC x_k + \bB \nabla g(\bC \tilde{x}_k) - \bB\bM z_f^*\\
    &= \bA x_k + \bB \underbrace{(\nabla g(\bC \tilde{x}_k) + \bM \bC \tilde{x}_k)}_{\nabla f(\bC x_k)}.
\end{align*}
This is the equation for the iterates $x_k$ of the algorithm defined by $(\bA,\bB,\bC)$ and $f$. Since $\tilde{x}_k$ goes to zero at rate $\rho$, so does $x_k$ go to $x_f^*$.
Finally, we argue that there cannot be an element of $S(\bM,\bL)$ with no critical point:
If there were an $f \in S(\bM,\bL)$ with two critical points, then the above arguments would prove convergence of algorithm \eqref{eq:algo} to both critical points, which cannot be true. Hence, there exists no such function in $S(\bM,\bL)$. Consequently, Theorem \ref{thm:ExMinCoCa}  guarantees that any function in $S(\bM,\bL)$ has a critical point  and thus $(\bA,\bB,\bC)$ solve Problem \ref{problem1}. (At this point the forward reference to Theorem \ref{thm:ExMinCoCa} can only be avoided by considerable effort. Also note that the proof of Theorem \ref{thm:ExMinCoCa} does in no way require Theorem \ref{thm:ConvTrInv}.) \qed

% As $f \in S(\bM,\bL)$ was chosen arbitrarily, this result holds for all $f \in S(\bM,\bL)$ and hence, $(\bA,\bB,\bC)$ solve Problem \ref{problem1}. \qed

\subsection{Proof of Theorem \ref{lem:LipschitzConv}}

 2) $\Rightarrow$ 1):
	The key to prove this statement is that the second term in inequality 1) can be written as the following integral:
	\begin{align*}
	    &\int_0^1 (\nabla f(z_1 + \tau (z_2-z_1)) - \nabla f(z_1))^T(z_2-z_1)\diff \tau\\
	    &= f(z_2)-f(z_1) - (\nabla f(z_1))^T(z_2-z_1).
	\end{align*}
	Using 2), the integrand can be upper and lower bounded as follows
	\begin{align*}
	    0&\leq (\nabla f(z_1 + \tau (z_2-z_1)) - \nabla f(z_1))^T(z_2-z_1)\\ &\leq \frac{1}{\tau} \| \tau (z_1-z_2)\|_{\bL}^2 = \tau \| z_1-z_2\|_{\bL}^2
	\end{align*}
	which implies 1).
	
 1) $\Rightarrow$ 3):
	Let $f\in C^1(\bbR^d)$ fulfil 1). Define $g(z) = f(z) - (\nabla f(z_1))^Tz$. Then $g \in S(0,\bL)$ and $\nabla g(z_1) = 0$. Thus $z_1$ is a minimizer of $g$ and we have
	\begin{align*}
	    g(z_1) - g(z_2) &\leq g(z_2 - \bA \nabla g(z_2)) - g(z_2)\\ &\overset{1)}{\leq}  \frac{1}{2}\|\bA \nabla g(z_2)\|_{\bL}^2 - \nabla g(z_2)^T \bA \nabla g(z_2),
	\end{align*}
	for any matrix $\bA \in \bbR^{d\times d}$ or equivalently 
	\begin{align*}
		\nabla g(z_2)^T \bA \nabla g(z_2) - \frac{1}{2}\|\bA \nabla g(z_2)\|_{\bL}^2 \leq g(z_2) - g(z_1).
	\end{align*}
	Now, we substitute  $g(z_2) = f(z_2) - (\nabla f(z_1))^Tz_2$:
	\begin{align*}
		 (\nabla f(z_1) - \nabla f(z_2))^T & \bA (\nabla f(z_1) - \nabla f(z_2))\\
		&- \frac{1}{2}\|\bA (\nabla f(z_1) - \nabla f(z_2))\|_{\bL}^2\\
		 \leq f(z_2) &- f(z_1) + (\nabla f(z_1))^T(z_1 - z_2).
	\end{align*}
	For $\bA = \bL^\dagger$, this is equivalent to
	\begin{align*}
		\frac{1}{2}\|\nabla f(z_1) &- \nabla f(z_2)\|_{\bL^\dagger}^2\\
		& \leq f(z_2) - f(z_1) + (\nabla f(z_1))^T(z_1 - z_2).
	\end{align*}
	In the case $\bA = r\bPi_{\ker \bL}$, the result is
	\begin{align*}
		r(\nabla f(z_1) - & \nabla f(z_2))^T \bPi_{\ker \bL} (\nabla f(z_1) - \nabla f(z_2))\\
		&\leq f(z_2) - f(z_1) + (\nabla f(z_1))^T(z_1 - z_2),
	\end{align*}
	which implies $\bPi_{\ker \bL} (\nabla f(z_1) - \nabla f(z_2)) = 0$, because $r$ can be chosen arbitrarily large. 
	
 3) $\Rightarrow$ 4):
	Adding the following inequalities %\ce{do not switch font size}
	$
		\frac{1}{2}\|\nabla f(z_1) - \nabla f(z_2)\|_{\bL^\dagger}^2
		\leq f(z_2) - f(z_1) + (\nabla f(z_1))^T(z_1 - z_2)
	$
	and
	$
		\frac{1}{2}\|\nabla f(z_1) - \nabla f(z_2)\|_{\bL^\dagger}^2
		\leq f(z_1) - f(z_2) + (\nabla f(z_2))^T(z_2 - z_1)
	$
	yields inequality in 4).
	
 4) $\Rightarrow$ 2):
	Let $f \in C^1(\bbR^d)$ fulfil 4). Then
	\begin{align*}
    	\sqrt{\bL} \sqrt{\bL^{\dagger}} (\nabla f(z_1) - \nabla f(z_2)) = (\nabla f(z_1) - \nabla f(z_2))
	\end{align*}
	holds for all $z_1,z_2 \in \bbR^d$, because $\bPi_{\ker \bL} (\nabla f(z_1) - \nabla f(z_2)) = 0$ implies, that $\nabla f(z_1) - \nabla f(z_2)$ is in the image of $\bL$. This observation can be used to derive the bound using the Cauchy-Schwarz-Inequality (CSI)
	{\small
	\begin{align*}
	    \| \nabla f(z_1) - \nabla f(z_2)\|_{\bL^\dagger}^2 \overset{4)}{\leq} (\nabla f(z_1) - \nabla f(z_2))^T (z_1 - z_2)&\\
		= (\nabla f(z_1) - \nabla f(z_2))^T \sqrt{\bL^{\dagger}} \sqrt{\bL} (z_1 - z_2)&\\
		\overset{\text{CSI}}{\leq} \| \nabla f(z_1) - \nabla f(z_2)\|_{\bL^\dagger} \| z_1 - z_2\|_{\bL},&
	\end{align*}}
	which implies $\| \nabla f(z_1) - \nabla f(z_2)\|_{\bL^\dagger} \leq \| z_1 - z_2\|_{\bL}$. Now, $f$ fulfils 2), because
	\begin{align*}
	    (\nabla f(z_1) - \nabla & f(z_2))^T (z_1 - z_2)\\ &\overset{\mathrm{CSI}}{\leq} \| \nabla f(z_1) - \nabla f(z_2)\|_{\bL^\dagger} \| z_1 - z_2\|_{\bL}\\
		    & \leq \| z_1 - z_2\|_{\bL}^2. \hspace{38mm} \qed
	\end{align*}

\subsection{Proof of Lemma \ref{lem:4.1EigCond}}

1) $\Rightarrow$ 2):
	Let $\bQ$ with $\bM \preceq \bQ \preceq \bL$ be given and let $(\lambda^{(\bM)}_i)_{i=1}^d$, $(\lambda^{(\bQ)}_i)_{i=1}^d$, $(\lambda^{(\bL)}_i)_{i=1}^d$ be the eigenvalues of those matrices in ascending order.
	It follows from $\bM \preceq \bQ \preceq \bL$ and the theorem of Courant-Fischer that
	\begin{align*}
		\lambda^{(\bM)}_1 \leq \lambda^{(\bQ)}_1 \leq \lambda^{(\bL)}_1, \ldots ,\lambda^{(\bM)}_d\leq \lambda^{(\bQ)}_d\leq \lambda^{(\bL)}_d
	\end{align*}
	holds. Since $\lambda^{(\bM)}_i$ and $\lambda^{(\bL)}_i$ always have the same sign and are not equal to zero by assumption, the values $\lambda^{(\bQ)}_i$ cannot be zero for any $i$. Hence, no eigenvalue of $\bQ$ can be zero and hence, $\bQ$ is invertible.

2) $\Rightarrow$ 3):
	To show the first statement, consider the case $\bQ = \frac{1}{2}(\bM + \bL)$. Then, it holds that $\bM \preceq \bQ \preceq \bL$ and hence, $\bQ = \frac{1}{2}(\bM + \bL)$ is invertible.
	To show the second statement, consider the case $\bQ = \frac{1}{2}(\bM + \bL) + \frac{\alpha}{2}(\bL- \bM)$. 
	For $\alpha \in [-1,1]$, it holds that $\bM \preceq \bQ \preceq \bL$ and thus
	\begin{align*}
		0&\neq \det \left( \frac{1}{2}(\bM + \bL) + \frac{\alpha}{2}(\bL- \bM) \right) & \forall \alpha \in [-1,1].
	\end{align*}
	By non-singularity of $(\bM + \bL)$, the factor $\det \frac{1}{2}(\bM + \bL)$ can be pulled out of the above expression, which gives
	\begin{align*}
		0&\neq \det \left(\frac{1}{2}(\bM + \bL)\right) \det \left( \bI + \alpha (\bM + \bL)^{-1} (\bL- \bM) \right)
	\end{align*}
	and consequently
	\begin{align*}
		0&\neq \det \left( \bI + \alpha (\bM + \bL)^{-1} (\bL- \bM) \right) & \forall \alpha \in [-1,1].
	\end{align*}
	This implies, that $(\bM + \bL)^{-1} (\bL- \bM)$ cannot have an eigenvalue in $\bbR\setminus ]-1,1[$. However, because $(\bM + \bL)^{-1} (\bL- \bM)$ is similar to the symmetric matrix $\sqrt{\bL- \bM} (\bM + \bL)^{-1} \sqrt{\bL- \bM}$, all of its eigenvalues have to be real. (Note that $\sqrt{\bL- \bM}$ exists because $\bL- \bM$ is positive semi-definite.) Hence, all eigenvalues of $(\bM + \bL)^{-1} (\bL- \bM)$ have to be in $]-1,1[$ and thus also $\rho ((\bM + \bL)^{-1} (\bL- \bM)) < 1$ holds.

3) $\Rightarrow$ 4):
	Suppose that $\bM$ is not invertible, i.e. there exists a vector $z \in \bbR^d\setminus \{0\}$ with $\bM z = 0$. Then
	\begin{align*}
		(\bL + \bM)z = (\bL- \bM) z
		\Rightarrow z = (\bL + \bM)^{-1}(\bL- \bM) z
	\end{align*}
	implies that $z$ is an eigenvector to the eigenvalue 1 of $(\bL + \bM)^{-1}(\bL- \bM)$, which contradicts $\rho ((\bM + \bL)^{-1} (\bL- \bM)) < 1$. Hence $\bM$ is non-singular.
	
	Next we show $\sigma (\bM^{-1} \bL) \subsEq \bbR_{> 0}$. Consider the identity
	\begin{align*}
		(\bL + \bM)^{-1}(\bL - \bM) &= \bI - 2(\bL + \bM)^{-1}\bM\\
		&= \bI - 2(\bM^{-1}\bL + \bI)^{-1}.
	\end{align*}
	Suppose, that $\bM^{-1}\bL$ has an eigenvalue $\lambda$ with associated eigenvector $v$. Then $\bM^{-1}\bL + \bI$ has eigenvalue $\lambda + 1$ with eigenvector $v$ and $(\bM^{-1} \bL + \bI)^{-1}$ has eigenvalue $\frac{1}{\lambda + 1}$ with eigenvector $v$. Thus
	\begin{align*}
		(\bL + \bM)^{-1}(\bL - \bM) v &=
		(\bI - 2(\bM^{-1}\bL + \bI)^{-1})v\\
		&= v - \frac{2}{\lambda + 1} v= \frac{\lambda -1}{\lambda + 1}v.
	\end{align*}
	Hence, $\frac{\lambda -1}{\lambda + 1}$ is an eigenvalue of $(\bL + \bM)^{-1}(\bL - \bM)$ and thus it is in $]-1,1[$. This implies $\lambda \in \bbR_{>0}$. Hence $\sigma (\bM^{-1}\bL) \subsEq \bbR_{> 0}$ holds true.

 4) $\Rightarrow$ 5):
	Suppose, that $\bL\bM^{-1}$ has only positive eigenvalues. Then there exists a symmetric positive definite matrix $\bP \in \bbR^{d\times d}$ such that the Lyapunov inequality
	\begin{align*}
		\bP \bL \bM^{-1} + \bM^{-1}\bL\bP \succ 0
	\end{align*}
	is satisfied. A congruence transform with $\bM$ gives
	\begin{align*}
	    \bM \bP \bL  + \bL\bP\bM \succ 0.
	\end{align*}
	By Lemma \ref{lem:congruence} we can infer that $\bM$ and $\bL$ are congruent.

5) $\Rightarrow$ 1):
	By Sylvester's Law of Inertia, matrices have the same eigenvalue signature, if and only if they are congruent. \qed

\subsection{Proof of Proposition \ref{prop:4.1GradContraction}}

We prove the contraction property of the map $\phi:z \mapsto z - 2(\bM + \bL)^{-1}\nabla f(z)$, by using the norm $\| z\|_{\bP}^2 = z^T\bP z$, where $\bP = (\bL + \bM)((\bL - \bM)^{\dagger} + r\bPi_{\ker (\bL-\bM)})(\bL + \bM)$.
In a first step, rewrite $\phi$ as:
\begin{align*}
	\phi(z) &= (\bL + \bM)^{-1}(\bL + \bM) z - 2(\bL + \bM)^{-1}\nabla f(z)\\
	&= (\bL + \bM)^{-1}\left((\bL - \bM) z - 2(\nabla f(z) - \bM z)\right)\\
	&= (\bL + \bM)^{-1}\left((\bL - \bM) z - 2\nabla g(z)\right)
\end{align*}
with $g \in S(0,\bL-\bM)$ defined by $g(z) := f(z) - \frac{1}{2}z^T\bM z$. 
Consider now
{\small
\begin{align*}
	&\|\phi(z_1) - \phi(z_2)\|^2_{\bP} \\
	&=
	\|((\bL - \bM)^{\dagger} + r\bPi_{\ker (\bL-\bM)})^{\frac{1}{2}}(\bL + \bM)(\bL + \bM)^{-1}\\
	&\hspace{3.7mm}\left(
	(\bL - \bM) (z_1 - z_2) - 2(\nabla g(z_1) - \nabla g(z_2))\right)\|^2\\
	&\resizebox{\linewidth}{!}{$\overset{(\star)}{=}
	\| \sqrt{(\bL - \bM)^{\dagger}}\left(
	(\bL - \bM) (z_1 - z_2) - 2(\nabla g(z_1) - \nabla g(z_2))\right)\|^2$}\\
	&=
	\|\sqrt{\bL - \bM} (z_1 - z_2) - 2\sqrt{(\bL - \bM)^{\dagger}}(\nabla g(z_1) - \nabla g(z_2))\|^2\\
	&= \underbrace{\resizebox{0.95\linewidth}{!}{$4 \|\nabla g(z_1) - \nabla g(z_2)\|_{(\bL - \bM)^{\dagger}}^2 - 4(\nabla g(z_1) - \nabla g(z_2))^T(z_1 - z_2)$}}_{\le 0 }\\
	&\hspace{3.7mm} + \|\sqrt{\bL - \bM}(z_1 - z_2)\|^2 \leq \|z_1 - z_2\|^2_{(\bL-\bM)}.
\end{align*}}
Concerning $(\star)$ notice, that the kernel projector has no contribution, 
since the products are all zero and
the under-braced expression being non-positive follows from Lemma \ref{lem:LipschitzConv}.
Finally, by Lemma \ref{lem:helpInequality} we know that for any $\rho > \rho_\mathrm{grad}$ there exists some $r \in \bbR_{>0}$ such that $\bL-\bM\preceq \rho^2 \bP$ holds. Hence, we can overestimate $\|z_1 - z_2\|^2_{(\bL-\bM)}$ by $\rho^2 \|z_1 - z_2\|^2_{\bP}$ (by choosing a sufficient value of $r$) resulting in the final estimate:
% Finally, by Lemma \ref{lem:helpInequality}, we can overestimate $\|z_1 - z_2\|^2_{(\bL-\bM)}$ as follows for any $\rho > \rho_\mathrm{grad}$ for some $r \in \bbR_{>0}$:
\begin{align*}
	\|\phi(z_1) - \phi(z_2)\|^2_{\bP} &\leq
	\|z_1 - z_2\|^2_{(\bL-\bM)} \leq \rho^2 \|z_1 - z_2\|^2_{\bP}.
\end{align*}
\qed

\subsection{Proof of Theorem \ref{thm:ExMinCoCa}}

Non-emptiness of $S(\bM,\bL)$ is equivalent to $\bM\preceq \bL$. It remains to show that the three statements in the theorem are equivalent under the condition $\bM \preceq \bL$.

\begin{itemize}
    \item 1) $\Rightarrow$ 2) and 1) $\Rightarrow$ 3):\\
	Assume $\bM \preceq_c \bL$ are non-singular. 
	Let $f \in S(\bM,\bL)$ be given. Then, by Proposition \ref{prop:4.1GradContraction}, the mapping
	\begin{align*}
		\phi : z \mapsto z - 2(\bM + \bL)^{-1}\nabla f(z)
	\end{align*}
	is a contraction on $\bbR^d$ and $(\bM + \bL)$ is non-singular. By 
	the Banach fixed point theorem the mapping $\phi$ has exactly one fixed point $z_f^*$ with
		$\phi(z_f^*) = z_f^* \Leftrightarrow \nabla f(z_f^*) = 0$.
	This implies 2) and 3).
	\item $\neg$1) $\Rightarrow$ $\neg$2) and $\neg$1) $\Rightarrow$ $\neg$3):\\
	Suppose that $\bM \preceq_c \bL$ does not hold or that either $\bM$ or $\bL$ or both are singular, but $\bM \preceq \bL$ holds (such that $S(\bM,\bL)$ is not empty). Then there exists $\bQ = \bQ^T \in \bbR^{d\times d}$ with $\bM \preceq \bQ \preceq \bL$ and $\det \bQ = 0$ by Lemma \ref{lem:4.1EigCond}.
	Let $v \in \bbR^d \setminus \{0\}$ be an element of the kernel of $\bQ$. Then the function $f_1\in S(\bM, \bL)$ defined by $f_1(z) = \frac{1}{2}z^T\bQ z + v^T z$ has no critical point with $\nabla f(z) = 0$, because otherwise
	\begin{align*}
		v^T \nabla f_1(z) = v^T(\bQ z + v) = \underbrace{v^T \bQ z}_{=0} +v^T v = \|v\|^2
	\end{align*}
	would have to be zero. At the same time, the function $f_2\in S(\bM, \bL)$ defined by $f_2(z) = \frac{1}{2}z^T\bQ z$ has infinitely many critical points with $\nabla f_2(z) = 0$, because any point $z = rv$ with $r \in \bbR$ is a critical point of $f$ by $\nabla f_2(z) = r\bQ v = 0$. \qed
\end{itemize}

\subsection{Proof of Lemma \ref{lem:4.4QuadBoundVf}}

 Step 1 (lower bound).
	The term $f(\bC x) - f(0) - \frac{1}{2}\nabla f(\bC x)^T \widetilde{\bL}^\dagger \nabla f(\bC x)$ can be lower bounded by the estimate
	{\small
	\begin{align*}
	0 &\leq f(\bC x) - f(0) - \frac{1}{2}\| \nabla f(\bC x) \|^2_{\widetilde{\bL}^{\dagger}} - (\nabla f(0))^T\bC x\\
	&= f(\bC x) - f(0) - \frac{1}{2}\nabla f(\bC x)^T \widetilde{\bL}^\dagger\nabla f(\bC x),
	\end{align*}}
	where the inequality sign follows from Lemma \ref{lem:LipschitzConv} and the equality sign follows from the fact $\nabla f(0) = 0$.
	This allows now the following lower bound on $V_f$:
	\begin{align*}
	V_f(x) &= \begin{pmatrix}
	x\\
	\nabla f(\bC x)
	\end{pmatrix}^T
	\begin{pmatrix}
	\bP_{11} & \bP_{12}\\
	\bP_{21} & \bP_{22}
	\end{pmatrix}
	\begin{pmatrix}
	x\\
	\nabla f(\bC x)
	\end{pmatrix}\\
	+& \underbrace{f(\bC x) - f(0) - \frac{1}{2}\nabla f(\bC x)^T \widetilde{\bL}^\dagger\nabla f(\bC x)}_{\geq 0}\\
	&\geq 
	\begin{pmatrix}
	x\\
	\nabla f(\bC x)
	\end{pmatrix}^T
	\begin{pmatrix}
	\bP_{11} & \bP_{12}\\
	\bP_{21} & \bP_{22}
	\end{pmatrix}
	\begin{pmatrix}
	x\\
	\nabla f(\bC x)
	\end{pmatrix}\\
	&\geq \lambda_\mathrm{min}(\bP) \left\| \begin{pmatrix}
	x\\
	\nabla f(\bC x)
	\end{pmatrix}\right\|^2\\
	&\geq \underbrace{\lambda_\mathrm{min}(\bP)}_{=:\alpha_f} \|x\|^2.
	\end{align*}
	
Step 2 (upper bound).
	The term $f(\bC x) - f(0) - \frac{1}{2}\nabla f(\bC x)^T \widetilde{\bL}^\dagger \nabla f(\bC x)$ can be upper bounded by the following estimates:
	\begin{align*}
	&f(\bC x) - f(0) \underbrace{- \frac{1}{2}\nabla f(\bC x)^T \widetilde{\bL}^\dagger\nabla f(\bC x)}_{\leq 0}\\
	&\overset{(\star)}{\leq} \underbrace{f(\bC x) - f(0) - (\nabla f(0))^T(\bC x - 0)}_{\overset{\text{Lemma } \ref{lem:LipschitzConv}}{\leq} \frac{1}{2}\|\bC x-0\|^2_{\widetilde{\bL}}}\\
	& \leq \frac{\|\bL-\bM\|}{2} \|\bC x - 0\|^2
	\end{align*}
	Note, that in $(\star)$ the term $(\nabla f(0))^T(\bC x - 0)$ can be added because $\nabla f(0) = 0$. This allows now the following upper bound on $V_f$:
	\begin{align*}
			V_f(x) =& \begin{pmatrix}
			x\\
			\nabla f(\bC x)
			\end{pmatrix}^T
			\begin{pmatrix}
			\bP_{11} & \bP_{12}\\
			\bP_{21} & \bP_{22}
			\end{pmatrix}
			\begin{pmatrix}
			x\\
			\nabla f(\bC x)
			\end{pmatrix}\\
			& + f(\bC x) - f(0) - \frac{1}{2}\nabla f(\bC x)^T \widetilde{\bL}^{\dagger}\nabla f(\bC x)\\
			\leq & \lambda_\mathrm{max}(\bP) \left\| \begin{pmatrix}
			x\\
			\nabla f(\bC x)
			\end{pmatrix}\right\|^2 + \frac{\|\widetilde{\bL}\|}{2} \|\bC x\|^2\\
			=& \lambda_\mathrm{max}(\bP)(\|x\|^2 + \underbrace{\|\nabla f(\bC x)\|^2}_{\leq \|\widetilde{\bL}\|^2 \|\bC x\|^2}) + \frac{\|\widetilde{\bL}\|}{2} \|\bC x\|^2\\
			\leq &  
			\left(\lambda_\mathrm{max}(\bP)(1+\|\widetilde{\bL}\|^2\|\bC\|^2) + \frac{\|\widetilde{\bL}\|\|\bC\|^2}{2}\right)
			\|x\|^2\\
			=& \beta_f \|x\|^2. \hspace{54mm} \qed
	\end{align*}

\subsection{Proof of Lemma \ref{lem:4.4EstDeriv}}

We define the abbreviations $w = \nabla f(\bC x)$, $w^+ = \nabla f(\bC x^+)$, $x^+ = \bA x + \bB w$ and $\widetilde{\bL} = \bL -\bM$. %Since the Lyapunov function candidates $V_f$ are dedicated for Problem \ref{problem2}, we have $f \in S(0,\bL- \bM)$ and thus by Lemma \ref{lem:LipschitzConv} $\bPi_{\bL - \bM} \nabla f(\bC x) = \bPi_{\bL - \bM} w = 0$. Consequently, we can add the term $r w^T \bPi_{\bL - \bM} w = 0$ for arbitrary $r \in \bbR$ to our Lyapunov function candidates resulting in
% \begin{align*}
%     V_f(x) =& \begin{pmatrix}
% 		x\\
% 		w
% 	\end{pmatrix}^T
% 	\begin{pmatrix}
% 		\bP_{11} & \bP_{12}\\
% 		\bP_{21} & \bP_{22}
% 	\end{pmatrix}
% 	\begin{pmatrix}
% 		x\\
% 		w
% 	\end{pmatrix}
% 	+ f(\bC x) - f(0)\\
% 	&- \frac{1}{2}\underbrace{w^T\left( (\bL - \bM)^{\dagger} + r \bPi_{\ker \bL - \bM} \right)w}_{=\frac{1}{2}\| w \|^2_{\bQ(r)}}.
% \end{align*}
With that the $\rho$-weighted increment of the Lyapunov function is
{\scriptsize
\begin{align*}
	&V_f(x^+) - \rho^2 V_f(x) =\\
%     &=
%     \begin{pmatrix}
% 	    x^+\\
% 	    w^+
% 	\end{pmatrix}^T
% 	\begin{pmatrix}
% 	    \bP_{11} & \bP_{12}\\
%     	\bP_{21} & \bP_{22}
% 	\end{pmatrix}
% 	\begin{pmatrix}
%     	x^+\\
% 	    w^+
% 	\end{pmatrix} + f(\bC x^+) - f(0) - \frac{1}{2}\| w^+ \|^2_{\bQ(r)}\\
% 	&-\rho^2\left(\begin{pmatrix}
% 	    x\\
% 	    w
% 	\end{pmatrix}^T
% 	\begin{pmatrix}
% 	    \bP_{11} & \bP_{12}\\
% 	    \bP_{21} & \bP_{22}
% 	\end{pmatrix}
% 	\begin{pmatrix}
%     	x\\
% 	    w
% 	\end{pmatrix}
% 	+ f(\bC x) - f(0) - \frac{1}{2}\| w \|^2_{\bQ(r)}\right)
% 	\\
	&=
	\begin{pmatrix}
		x^+\\
		w^+
	\end{pmatrix}^T
	\begin{pmatrix}
		\bP_{11} & \bP_{12}\\
		\bP_{21} & \bP_{22}
	\end{pmatrix}
	\begin{pmatrix}
		x^+\\
		w^+
	\end{pmatrix} 
	-\rho^2\begin{pmatrix}
		x\\
		w
	\end{pmatrix}^T
	\begin{pmatrix}
		\bP_{11} & \bP_{12}\\
		\bP_{21} & \bP_{22}
	\end{pmatrix}
	\begin{pmatrix}
		x\\
		w
	\end{pmatrix}\\
	&+ \underbrace{f(\bC x^+) - f(0) - \frac{1}{2}\| w^+ \|^2_{\widetilde{\bL}^\dagger}
		-\rho^2\left(f(\bC x) - f(0) - \frac{1}{2}\| w \|^2_{\widetilde{\bL}^\dagger}\right)}_{I}.
\end{align*}}
To upper bound expression $I$, we use the estimate
\begin{align*}
    &\underbrace{-\rho^2}_{\leq - \lambda}\underbrace{\left(f(\bC x) - f(0) - \frac{1}{2}\| w \|^2_{\widetilde{\bL}^\dagger}\right)}_{\geq 0}\\
    &\leq - \lambda \left(f(\bC x) - f(0) - \frac{1}{2}\| w \|^2_{\widetilde{\bL}^\dagger}\right),
\end{align*}
which we can use to obtain
\begin{align*}
	I &\leq (1-\lambda)\underbrace{\left(f(\bC x^+) - f(0) + \frac{1}{2}\| w^+ \|^2_{\widetilde{\bL}^\dagger}\right)}_{\overset{\text{Lemma } \ref{lem:LipschitzConv}}{\leq} (w^+)^T(\bC x^+ - 0)} \\
	&\hspace{3.7mm} + \lambda\underbrace{\left(f(\bC x^+) - f(\bC x) + \frac{1}{2}\| w^+ - w \|^2_{\widetilde{\bL}^\dagger}\right)}_{\overset{\text{Lemma } \ref{lem:LipschitzConv}}{\leq} (w^+)^T(\bC x^+ - \bC x)}\\
	&\hspace{3.7mm}-\underbrace{\frac{(2 - \lambda)}{2}\| w^+ \|^2_{\widetilde{\bL}^\dagger} - \frac{\lambda}{2}\| w^+ - w \|^2_{\widetilde{\bL}^\dagger} + \frac{\lambda}{2}\| w \|^2_{\widetilde{\bL}^\dagger}}_{ = (w^+)^T\widetilde{\bL}^\dagger(w^+ - \lambda w)}\\
	&\leq (1-\lambda)(w^+)^T\bC x^+ + \lambda(w^+)^T(\bC x^+ - \bC x)\\ &\hspace{3.7mm} -(w^+)^T\widetilde{\bL}^\dagger(w^+ - \lambda w)\\
	&= (w^+)^T\left(\bC x^+-\lambda \bC x - \widetilde{\bL}^\dagger(w^+ - \lambda w)\right).\\
%	&= (w^+)^T\left(\bC x^+-\lambda \bC x - (\bQ(r)w^+ - \lambda \bQ(0) w)\right).
\end{align*}
%Notice in the above calculation, that for all these estimates $\bQ w  = (\bL-\bM)^\dagger w$. Hence, Lemma  \ref{lem:LipschitzConv} can be applied with $\bQ$ instead of $(\bL-\bM)^\dagger$.\\
Now, this estimate for expression $I$ can be used to upper bound $V_f(x^+) - \rho^2 V_f(x)$ as follows:
\begin{align*}
	V(x^+) - \rho^2 V(x) &\leq \begin{pmatrix}
    	x^+\\
	    w^+
	\end{pmatrix}^T
	\begin{pmatrix}
	    \bP_{11} & \bP_{12}\\
	    \bP_{21} & \bP_{22}
	\end{pmatrix}
	\begin{pmatrix}
	    x^+\\
	    w^+
	\end{pmatrix}\\
	&\hspace{3.7mm} -\rho^2\begin{pmatrix}
	    x\\
	    w
	\end{pmatrix}^T
	\begin{pmatrix}
	    \bP_{11} & \bP_{12}\\
	    \bP_{21} & \bP_{22}
	\end{pmatrix}
	\begin{pmatrix}
	    x\\
	    w
	\end{pmatrix}\\
	+ (w^+)^T & \left(\bC x^+-\lambda \bC x - \widetilde{\bL}^\dagger( w^+ - \lambda w)\right),
\end{align*}
which corresponds to the inequality in Lemma \ref{lem:4.4EstDeriv}.
 \qed

\subsection{Proof of Theorem \ref{thm:4.5}}

First remember that Theorem \ref{thm:ConvTrInv} shows that an algorithm with parameters $(\bA,\bB,\bC)$ has convergence rate $\rho$ for $S(\bM,\bL)$ if 
an algorithm with parameters $(\widetilde{\bA},\bB,\bC)$, which satisfy the constraint \eqref{eq:Constraint}, has convergence rate $\rho$ for $S_0(0,\widetilde{\bL})$.
Hence, in the following we show convergence for $(\widetilde{\bA},\bB,\bC)$ and $S_0(0,\widetilde{\bL})$.
By Theorem \ref{thm:2.1LyapOpt}, an algorithm defined by $(\widetilde{\bA},\bB,\bC)$ is asymptotically stable and has convergence rate $\rho$, if there exists a Lyapunov function $V_f:\bbR^n\to\bbR$, such that 
\begin{align*}
	&\alpha_f \| x - x_f^*\|^2 \leq V_f(x) \leq \beta_f \| x - x_f^*\|^2, \\
	&V_f(x^+) - \rho^2 V_f(x)\leq 0 
\end{align*}
holds for all $x \in \bbR^n$ and $f \in S_0(0,\widetilde{\bL})$ with $\beta_f \geq \alpha_f > 0$.
The considered class of Lyapunov function candidates fulfills these requirements by Lemma \ref{lem:4.4QuadBoundVf} and Lemma \ref{lem:4.4EstDeriv} if
{\small
\begin{align*}
	&\left(
	\begin{array}{c}
	x\\
	w\\
	\hline
	x^+\\
	w^+
	\end{array}\right)^T
	\left(
	\begin{array}{cc|cc}
	-\rho^2 \bP_{11} & -\rho^2 \bP_{12} & 0 & 0\\
	-\rho^2 \bP_{21} & -\rho^2 \bP_{22} & 0 & 0\\
	\hline
	0 & 0 & \bP_{11} & \bP_{12}\\
	0 & 0 & \bP_{21} & \bP_{22}
	\end{array}
	\right)
	\left(
	\begin{array}{c}
	x\\
	w\\
	\hline
	x^+\\
	w^+
	\end{array}\right)\\
	&+
	\left(
	\begin{array}{c}
	x\\
	w\\
	\hline
	x^+\\
	w^+
	\end{array}\right)
	\left(
	\begin{array}{cc|cc}
	0 & 0 & 0 & -\frac{\lambda}{2}\bC^T\\
	0 & 0 & 0 & \frac{\lambda}{2}\widetilde{\bL}^\dagger\\
	\hline
	0 & 0 & 0 & \frac{1}{2}\bC^T\\
	-\frac{\lambda}{2}\bC & \frac{\lambda}{2}\widetilde{\bL}^\dagger & \frac{1}{2}\bC & -\widetilde{\bL}^\dagger
	\end{array}\right)
	\left(
	\begin{array}{c}
	x\\
	w\\
	\hline
	x^+\\
	w^+
	\end{array}\right)
\end{align*}}
is smaller than zero for all $x\in \bbR^n$, $w = \nabla f(\bC x)$, $w^+ = \nabla f(\bC x^+)$ and $x^+ = \bA x + \bB w$. At this point we can even improve the estimate by the observation that due to Lemma~\ref{lem:LipschitzConv}
\begin{align*}
	0 &= \bPi_{\ker \widetilde{\bL}} \nabla f(\bC x) = \bPi_{\ker \widetilde{\bL}} w,\\ 
	0 &= \bPi_{\ker \widetilde{\bL}} \nabla f(\bC x^+) = \bPi_{\ker \widetilde{\bL}} w^+
\end{align*}
hold true. This implies, that the term
\begin{align*}
    \left(
	\begin{array}{c}
	x\\
	w\\
	\hline
	x^+\\
	w^+
	\end{array}\right)
	\left(
	\begin{array}{cc|cc}
	0 & 0 & 0 & 0\\
	0 & -r\bPi_{\ker \widetilde{\bL}} & 0 & 0\\
	\hline
	0 & 0 & 0 & 0\\
	0 & 0 & 0 & -r\bPi_{\widetilde{\bL}}
	\end{array}\right)
	\left(
	\begin{array}{c}
	x\\
	w\\
	\hline
	x^+\\
	w^+
	\end{array}\right)
\end{align*}
is zero for all $r\in \bbR$ and can hence be added (as an additional multiplier) to the estimate.
Since the quantities $x,w,x^+,w^+$ are given by
\begin{align*}
	\left(
	\begin{array}{c}
	x\\
	w\\
	\hline
	x^+\\
	w^+
	\end{array}\right)
	=
	\left(
	\begin{array}{ccc}
	\bI_n & 0 & 0\\
	0 & \bI_d & 0\\
	\hline
	\widetilde{\bA} & \bB & 0\\
	0 & 0 & \bI_d
	\end{array}\right)
	\begin{pmatrix}
	x \\
	w\\
	w^+
	\end{pmatrix}
\end{align*}
negativity of $V_f(x^+) - \rho^2 V_f(x)$ follows now from inequality \eqref{eq:4.5_1}. Hence, \eqref{eq:4.5_1} implies that the weighted increment of the Lyapunov function is negative definite and, as a consequence, that the algorithm defined by $(\widetilde{\bA},\bB,\bC)$ has convergence rate $\rho$ for $S_0(0,\bL-\bM)$. \qed

\subsection{Proof of Theorem \ref{thm:4.6SynIneq}}

We need to show  that the  matrix inequality \eqref{eq:SynIneq} in the transformed variables $\hat{\bA},\hat{\bB},\bC,\bP$ is equivalent to \eqref{eq:4.5_1}.
The proof of this theorem works in two steps. The first step is to apply the Schur complement to \eqref{eq:4.5_1}. The second (key) step is to define a linearizing change of variables.

Step 1. First, define $\bZ$ as follows
{\small
\begin{align*}
    \resizebox{\linewidth}{!}{$
	\left(
	\begin{array}{ccc}
	\bI_n & 0 & 0\\
	0 & \bI_d & 0\\
	\hline
	\widetilde{\bA} & \bB & 0\\
	0 & 0 & \bI_d
	\end{array}\right)^T\hspace{-1mm}
	\left(
	\begin{array}{cc|cc}
	0 & 0 & 0 & -\frac{\lambda}{2}\bC^T\\
	0 & -r\bPi & 0 & \frac{\lambda}{2}\widetilde{\bL}^\dagger\\
	\hline
	0 & 0 & 0 & \frac{1}{2}\bC^T\\
	-\frac{\lambda}{2}\bC & \frac{\lambda}{2}\widetilde{\bL}^\dagger & \frac{1}{2}\bC & -\widetilde{\bL}^\dagger-r\bPi
	\end{array}\right)
	(\star)$}\\
% 	\left(
% 	\begin{array}{ccc}
% 	\bI_n & 0 & 0\\
% 	0 & \bI_d & 0\\
% 	\hline
% 	\widetilde{\bA} & \bB & 0\\
% 	0 & 0 & \bI_d
% 	\end{array}\right)
	\resizebox{\linewidth}{!}{$=
	\begin{pmatrix}
	0 & 0 & \frac{1}{2}\widetilde{\bA}^T\bC^T  -\frac{\lambda}{2}\bC^T \\
	0 & -r\bPi & \frac{1}{2} \bB^T\bC^T  + \frac{\lambda}{2}\widetilde{\bL}^\dagger \\
	\frac{1}{2}\bC\widetilde{\bA}  -\frac{\lambda}{2}\bC & \frac{1}{2} \bC\bB  + \frac{\lambda}{2}\widetilde{\bL}^\dagger & -\widetilde{\bL}^\dagger-r\bPi
	\end{pmatrix}=:\bZ.$}
\end{align*}}
 With $\bZ$, \eqref{eq:4.5_1} becomes
{\small
\begin{align*}
	& 
	\resizebox{\linewidth}{!}{$
	\left(
	\begin{array}{ccc}
	\bI_n & 0 & 0\\
	0 & \bI_d & 0\\
	\hline
	\widetilde{\bA} & \bB & 0\\
	0 & 0 & \bI_d
	\end{array}\right)^T
	\left(
	\begin{array}{cc|cc}
	-\rho^2 \bP_{11} & -\rho^2 \bP_{12} & 0 & 0\\
	-\rho^2 \bP_{21} & -\rho^2 \bP_{22} & 0 & 0\\
	\hline
	0 & 0 & \bP_{11} & \bP_{12}\\
	0 & 0 & \bP_{21} & \bP_{22}
	\end{array}
	\right)
	(\star) + \bZ$}
	\\
% 	&=
% 	\begin{pmatrix}
% 	\begin{pmatrix}
% 	\bI_n & 0\\
% 	0 & \bI_d
% 	\end{pmatrix}
% 	& 0\\
% 	\bP
% 	\begin{pmatrix}
% 	\widetilde{\bA} & \bB\\
% 	0 & 0
% 	\end{pmatrix}
% 	& \bP\begin{pmatrix}
% 	0 \\ \bI_d
% 	\end{pmatrix}
% 	\end{pmatrix}^T
% 	\begin{pmatrix}
% 	-\rho^2 \bP & 0\\
% 	0 & \bP^{-1}
% 	\end{pmatrix}
% 	(\star)
% 	+ \bZ\\
	&=
	\begin{pmatrix}
	\bI_n & 0 & 0\\
	0 & \bI_d & 0\\
	\bP_{11}\widetilde{\bA} & \bP_{11}\bB & \bP_{12}\\
	\bP_{21}\widetilde{\bA} & \bP_{21}\bB & \bP_{22}
	\end{pmatrix}^T
	\begin{pmatrix}
	-\rho^2 \bP & 0\\
	0 & \bP^{-1}
	\end{pmatrix}
	(\star)
	+ \bZ\prec 0.
\end{align*}}
The matrix $\bP$ is positive definite, by assumption of Theorem \ref{thm:4.5} and as a consequence of the matrix inequality from Theorem \ref{thm:4.6SynIneq}.
Hence, this algebraic manipulation allows to apply the Schur complement, which states that the above inequality is equivalent to
\begin{align*}
    \resizebox{\linewidth}{!}{$
	\begin{pmatrix}
	-\rho^2 \bP_{11} & -\rho^2 \bP_{12} & * & * & *\\
	-\rho^2 \bP_{21} & -\rho^2 \bP_{22}-r\bPi & * & * & *\\
	\frac{1}{2}\bC\widetilde{\bA}  -\frac{\lambda}{2}\bC & \frac{1}{2} \bC\bB  + \frac{\lambda}{2}\widetilde{\bL}^\dagger & -\widetilde{\bL}^\dagger-r\bPi & * & *\\
	\bP_{11}\widetilde{\bA} & \bP_{11} \bB & \bP_{12} & -\bP_{11} & -\bP_{12}\\
	\bP_{21} \widetilde{\bA} & \bP_{21}\bB & \bP_{22} & -\bP_{21} & -\bP_{22}
	\end{pmatrix}$}
\end{align*}
being negative definite. %The requirement of the Schur complement Lemma, that $\bP$ must be positive definite, is fulfilled by the assumptions of Theorem \ref{thm:4.5} and implied by the matrix inequality in Theorem \ref{thm:4.6SynIneq}. Hence, it is always fulfilled.

Step 2. If we have a solution $(\hat{\bA},\hat{\bB},\ldots)$
of \eqref{eq:SynIneq} and constraint \eqref{eq:ConstraintSyn}, then we can just substitute $\widetilde{\bA} = \bP_{11}^{-1}\hat{\bA}, \bB = \bP_{11}^{-1}\hat{\bB}$ into \eqref{eq:SynIneq} and we see that we obtain
the above inequality and hence a solution of \eqref{eq:4.5_1}.
This solution also  satisfies constraint \eqref{eq:Constraint} since
\begin{align*}
    \bC(\widetilde{\bA} - \bI_n)^{-1} \bB \bM &= \bC(\widetilde{\bA} - \bI_n)^{-1} \bP_{11}^{-1}\hat{\bB} \bM\\
    &= \bC(\bP_{11}\widetilde{\bA} - \bP_{11})^{-1} \hat{\bB} \bM\\
    &= \bC(\hat{\bA} - \bP_{11})^{-1} \hat{\bB} \bM\\
    &= \bC \bJ_1^{T} = \bI_d.
\end{align*}
On the other hand, if we are given a solution of \eqref{eq:4.5_1}, \eqref{eq:Constraint} with $\bP \succ 0$ and we want to construct a solution of \eqref{eq:SynIneq} by substituting $\hat{\bB} = \bP_{11}\bB$, $\hat{\bA} = \bP_{11}\widetilde{\bA}$ and by expressing all the nonlinear expressions $\bC \widetilde{\bA}, \bP_{21} \widetilde{\bA}, \bP_{11}\widetilde{\bA}, \bC\bB, \bP_{21}\bB, \bP_{11}\bB$ in terms of $\hat{\bA}$ and $\hat{\bB}$, we cannot guarantee that \eqref{eq:ConstraintSyn} holds.
However, in the following we show that there exists a state transformation of the algorithm such that this can be indeed guaranteed. Hence, any solution of \eqref{eq:4.5_1}, \eqref{eq:Constraint} is a solution of \eqref{eq:SynIneq}, \eqref{eq:ConstraintSyn} by an appropriate coordinate transformation.

If there exists a transformation (non-singular) matrix $\bT$ such that the transformed variables $\widetilde{\bA}' = \bT^{-1}\widetilde{\bA}\bT$, $\bB' = \bT^{-1}\bB$, $\bC' = \bC\bT$, $\bP_{11}' = \bT^T\bP_{11}\bT$, $\bP_{12}' = \bT^T\bP_{12}$, $\bP_{21}' = \bP_{21}\bT$, $\bP_{22}' = \bP_{22}$ fulfill
\begin{align*}
	(\widetilde{\bA}' - \bI_n)\bJ_1^T = \bB'\bM,\hspace{4mm} \bJ_2\bP_{11}' = \bC', \hspace{4mm} \bJ_3\bP_{11}' = \bP_{21}',
\end{align*}
then we have
%all the blocks can be rewritten in terms of $\hat{\bA}'$ and $\hat{\bB}'$, as follows
\begin{align*}
	\begin{pmatrix}
	\bC'\widetilde{\bA}'\\
	\bP_{21}'\widetilde{\bA}'\\
	\bP_{11}'\widetilde{\bA}'
	\end{pmatrix} 
	=
	\begin{pmatrix}
	\bJ_2\\
	\bJ_3\\
	\bI_n
	\end{pmatrix}\hat{\bA}',
	\hspace{7mm}
	\begin{pmatrix}
	\bC'\bB'\\
	\bP_{21}'\bB'\\
	\bP_{11}'\bB'
	\end{pmatrix} 
	=
	\begin{pmatrix}
	\bJ_2\\
	\bJ_3\\
	\bI_n
	\end{pmatrix}\hat{\bB}'
\end{align*}
and the transformed variables still form a solution of inequality \eqref{eq:4.5_1}. The arguments from Step 1 show that in this case $\widetilde{\bA}'$, $\bB'$, $\bC'$ and $\bP'$ form also a solution of \eqref{eq:SynIneq} and by substituting $\widehat{\bA}'$ and $\widehat{\bB}'$ for the nonlinear terms it becomes clear that there exists a solution to \eqref{eq:SynIneq}, \eqref{eq:ConstraintSyn} from Theorem \ref{thm:4.6SynIneq}.\\
Such a transformation $\bT$ must now fulfil the constraints
\begin{align*}
	\bJ_2\underbrace{\bT^T\bP_{11}\bT}_{=\bP_{11}'} = \underbrace{\bC\bT}_{ = \bC'}&,\hspace{5mm} \bJ_3\bT^T\bP_{11}\bT = \underbrace{\bP_{21}\bT}_{=\bP_{21}'},\\ \underbrace{\bT^{-1}(\widetilde{\bA} - \bI_n)\bT}_{=\widetilde{\bA}'-\bI_n}\bJ_1^T &= \underbrace{\bT^{-1}\bB}_{=\bB'}\bM.
\end{align*}
Rearranging and canceling terms above gives
{\small
\begin{align*}
    \resizebox{\linewidth}{!}{$
	\bT\bJ_1^T = (\widetilde{\bA}-\bI_n)^{-1}\bB\bM, \hspace{3mm} \bJ_2\bT^T = \bC\bP_{11}^{-1},\hspace{3mm} \bJ_3\bT^T = \bP_{21}\bP_{11}^{-1}.$}
\end{align*}}
For the choice $\bJ_1 = (\bI_d~0), \bJ_2 = (0_d~\bI_d~0), \bJ_3 = (0_d~0_d~\bI_d~0)$, these equations have the solution
\begin{align*}
	\bT = \begin{pmatrix}
	(\widetilde{\bA}-\bI_n)^{-1}\bB\bM & \bP_{11}^{-T}\bC^{T} & \bP_{11}^{-T}\bP_{21}^T & \bT_4
	\end{pmatrix},
\end{align*}
provided, that $n \geq 3d$. %\ce{If we have enough time, we should think if there is a way to avoid $n \geq 3d$}
It remains to show that the transformation is non-singular.
Notice that  $(\widetilde{\bA}-\bI_n)^{-1}\bB\bM$ and $\bC$ must have full rank because of  $\bC(\widetilde{\bA}-\bI_n)^{-1}\bB\bM = \bI_d$. Moreover $\bP_{11},\bP_{21}$ can be slightly perturbed (without violating the strict definiteness of $\bP$ and the matrix inequality \eqref{eq:4.5_1}), such that $\left((\widetilde{\bA}-\bI_n)^{-1}\bB\bM ~~ \bP_{11}^{-T}\bC^{T} ~~ \bP_{11}^{-T}\bP_{21}^T\right)$ has full rank too. Finally, $\bT_4 \in \bbR^{n\times n-3d}$ can be chosen such that $\bT$ is non-singular. 
Hence, all constraints of Theorem \ref{thm:4.6SynIneq} are satisfied by construction of $\bT$, where $\bC \bJ_1 = \bI_d$ is implied by \eqref{eq:Constraint}.

Consequently, it is possible to construct solutions related to Theorem \ref{thm:4.6SynIneq} from solutions related to Theorem \ref{thm:4.5} and vice versa. \qed

\subsection{Proof of Theorem \ref{thm:Existence of solutions}}

Again, we introduce the abbreviations $\widetilde{\bL} = \bL - \bM$ and $\bPi = \bPi_{\ker \bL - \bM}$. In this proof, it will be necessary to find explicit solutions for the matrix inequality \eqref{eq:4.5_1} from Theorem \ref{thm:4.5}. Therefore, it is purposeful to multiply out the matrix products in this inequality for $\lambda = 0$, which gives:
\begin{align*}
    \resizebox{\linewidth}{!}{$
	\begin{pmatrix}
	\widetilde{\bA}^T\bP_{11}\widetilde{\bA} - \rho^2 \bP_{11} & \widetilde{\bA}^T\bP_{11} \bB - \rho^2 \bP_{12} & \widetilde{\bA}^T(\bP_{12} + \frac{1}{2}\bC^T)\\
	\bB \bP_{11}\widetilde{\bA} - \rho^2 \bP_{21} & \bB^T \bP_{11}\bB - \rho^2\bP_{22}-r\bPi & \bB^T (\bP_{12} + \frac{1}{2}\bC^T)\\
	(\bP_{21} + \frac{1}{2}\bC)\widetilde{\bA} & (\bP_{21} + \frac{1}{2}\bC)\bB & \bP_{22} - \widetilde{\bL}^\dagger -r\bPi
	\end{pmatrix}.$}
\end{align*}

i) $\Rightarrow$ ii): 
	This step will be quite lengthy. We will show, that the matrices $(\widetilde{\bA},\bB,\bC)$ given by
	\begin{align*}
	\widetilde{\bA} &= \bA + \bB\bM\bC = \bI_d - 2(\bL + \bM)^{-1}\bM\\
	&= (\bL+\bM)^{-1}(\bL - \bM)\\
	\bB &= -2(\bL+\bM)^{-1}\\
	\bC &= \bI_d
	\end{align*}
	fulfill all the convergence rate conditions of Theorem \ref{thm:4.5} for an arbitrary given $\rho \in ]\rho_\mathrm{grad},1[$.	Here, the matrix $\widetilde{\bA}$ fulfils the Lyapunov inequality
	\begin{align*}
		\widetilde{\bA}^T\widetilde{\bP} \widetilde{\bA} - \rho^2 \widetilde{\bP} \prec 0
	\end{align*}
	for $\widetilde{\bP} := (\bL + \bM)\left((\bL - \bM)^{\dagger}+r \bPi \right)(\bL+\bM)$ and large enough $r \in \bbR_{>0}$ by Lemma \ref{lem:helpInequality}, since $\widetilde{\bA}^T\widetilde{\bP} \widetilde{\bA} = \bL - \bM$.
	To show, that the convergence conditions from Theorem \ref{thm:4.5} are met we choose $\lambda = 0$ and the following value for $\bP$: 
		\begin{align*}
		\resizebox{\linewidth}{!}{$
		\begin{pmatrix}
		\bP_{11} & \bP_{12}\\
		\bP_{21} & \bP_{22}
		\end{pmatrix}
		=
		\begin{pmatrix}
		\frac{\rho^2}{4} \left(\widetilde{\bP} - \varepsilon  (\bL + \bM)^2\right) & -\frac{1}{2}\bI_d\\
		-\frac{1}{2}\bI_d & \widetilde{\bL}^\dagger+r\bPi - \frac{\varepsilon}{2} \bI_d
		\end{pmatrix},$}
		\end{align*}
		where $\varepsilon > 0$, and $r \in \bbR$ is the same as above.
		There are three things to show:
		\begin{enumerate}[1)]
			\item The constraint \eqref{eq:Constraint} of Theorem \ref{thm:4.5} is satisfied for $\widetilde{\bA},\bB,\bC$.
			\item For large enough $r$ and small enough $\varepsilon$, $\bP$ solves the matrix inequality \eqref{eq:4.5_1} of Theorem \ref{thm:4.5}.
			\item For large enough $r$ and small enough $\varepsilon$, $\bP$ is positive definite.
		\end{enumerate}
		Verifying 1) can be done by a simple calculation of formulas in the constraint.\\
		We will now show 2). Note that $(\bP_{21} + \frac{1}{2}\bC) = \frac{1}{2}(\bI_d - \bI_d) = 0$ holds, which is why \eqref{eq:4.5_1} from Theorem \ref{thm:4.5} simplifies to
		\begin{align*}
		\resizebox{\linewidth}{!}{$
		\begin{pmatrix}
		\widetilde{\bA}^T\bP_{11}\widetilde{\bA} - \rho^2 \bP_{11} & \widetilde{\bA}^T\bP_{11} \bB - \rho^2 \bP_{12} & 0\\
		\bB \bP_{11}\widetilde{\bA} - \rho^2 \bP_{21} & \bB^T \bP_{11}\bB - \rho^2\bP_{22} -r\bPi & 0\\
		0 & 0 & -\frac{\varepsilon}{2} \bI_d
		\end{pmatrix} \prec 0.$}
		\end{align*}
		Here it is left to show, that the left upper $2\times 2$ block can be made negative definite by choosing $r$ big and $\varepsilon$ small. This is done by dividing the matrix inequality by $\frac{\rho^2}{4}$ and calculating the entries of the left upper blocks:\\
		The first block is
		\begin{align*}
		&\frac{4}{\rho^2}\left(\widetilde{\bA}^T\bP_{11}\widetilde{\bA} - \rho^2 \bP_{11}\right)\\
		&=
		\widetilde{\bA}^T\left(\widetilde{\bP} - \frac{\varepsilon}{2}  (\bL + \bM)^2\right)\widetilde{\bA} - \rho^2 \left(\widetilde{\bP} - \varepsilon  (\bL + \bM)^2\right)\\
		&= 
		\widetilde{\bA}^T\widetilde{\bP}\widetilde{\bA} - \rho^2 \widetilde{\bP} - \varepsilon\left( (\bL-\bM)^2 - \rho^2 (\bL + \bM)^2\right).
		\end{align*}
		The second block is
		\begin{align*}
		&\frac{4}{\rho^2}\left(\widetilde{\bA}^T\bP_{11} \bB - \rho^2 \bP_{12}\right)\\
		&=
		2\bI_d + \widetilde{\bA}^T \left(\widetilde{\bP} - \varepsilon  (\bL + \bM)^2\right) \bB\\
		&= 2\bI_d -2(\bL - \bM) \left((\bL - \bM)^{\dagger}+r \bPi \right) - 2\varepsilon(\bL - \bM)\\
		&\resizebox{\linewidth}{!}{$= 
		2 \bI_d -2\underbrace{(\bL - \bM) (\bL - \bM)^{\dagger}}_{\overset{\eqref{eq:PseudoInvProd}}{=}\bPi_{\im (\bL-\bM)}} -2r \underbrace{(\bL- \bM)\bPi}_{ = 0} - 2\varepsilon(\bL - \bM)$}\\
		&= 2\bI_d -2\bPi_{\im (\bL-\bM)} - 2\varepsilon(\bL - \bM)\\
		&\overset{\eqref{eq:ProjektorSum}}{=} 2 \bPi_{\ker (\bL-\bM)} - 2\varepsilon(\bL - \bM).
		\end{align*}
		The third block is:
		\begin{align*}
		&\frac{4}{\rho^2}\left(\bB^T \bP_{11} \bB - \rho^2\bP_{22}- r\bPi\right)\\
		&= \bB^T \widetilde{\bP} \bB - 4\varepsilon\bI_d - 4\left(\widetilde{\bL}^\dagger + r\bPi - \frac{\varepsilon}{2} \bI_d\right) -\frac{4}{\rho^2} r\bPi\\
		&= 4(\widetilde{\bL}^\dagger + r\bPi) - 4\varepsilon\bI_d - 4\left(\widetilde{\bL}^\dagger+r\bPi - \frac{\varepsilon}{2} \bI_d\right) -\frac{4}{\rho^2} r\bPi\\
		&= 4 \left(\widetilde{\bL}^\dagger + r\bPi - \widetilde{\bL}^\dagger-r\bPi\right) - 2\varepsilon\bI_d -\frac{4}{\rho^2} r\bPi.
		\end{align*}
		The calculation of these blocks reveals, that the upper $2\times 2$ block is		
		\begin{align*}
		\resizebox{\linewidth}{!}{$
		\begin{pmatrix}
		\widetilde{\bA}^T\widetilde{\bP}\widetilde{\bA} - \rho^2 \widetilde{\bP} - \varepsilon\left( (\bL-\bM)^2 - \rho^2 (\bL + \bM)^2\right) & 2 \bPi-2\varepsilon (\bL-´\bM)\\
		2 \bPi-2\varepsilon (\bL-´\bM) & -2\varepsilon \bI_d -\frac{4}{\rho^2}r\bPi\\
		\end{pmatrix},$}
		\end{align*}
		which is negative definite for $\varepsilon > 0$ small enough and $r$ big enough.\\
		Now it is left to show 3), namely that $\bP$ is positive definite for small enough $\varepsilon$ and large enough $r$. Therefore, we can show that $\bP$ is positive definite for $\varepsilon = 0$. Then it will also be positive definite for the small perturbation with $\varepsilon > 0$. By the Schur complement, the matrix $\bP$ for $\varepsilon = 0$ is positive definite if and only if:
		\begin{align*}
		0&\prec \widetilde{\bL}^\dagger + r\bPi\\
		0 &\prec \frac{\rho^2}{4} \widetilde{\bP} - \left( -\frac{1}{2}\bI_d \right) \underbrace{(\widetilde{\bL}^\dagger + r\bPi)^{-1}}_{\overset{\eqref{eq:InverseWithPseudoInverse2}}{=} (\bL - \bM) + \frac{1}{r}\bPi} \left( -\frac{1}{2}\bI_d \right)\\
		&= \frac{\rho^2}{4} \widetilde{\bP} - \frac{1}{4}(\bL - \bM) - \frac{1}{4r}\bPi.
		\end{align*}
		Since $\rho > \rho_\mathrm{grad}$, the matrix $\rho^2 \widetilde{\bP} - (\bL - \bM)$ is positive definite by Lemma \ref{lem:helpInequality} and the matrix $\widetilde{\bL}^\dagger + r \bPi$ is positive definite by construction. Thus, $\frac{\rho^2}{4} \widetilde{\bP} - \frac{1}{4}(\bL - \bM) - \frac{1}{4r}\bPi$ is positive definite for large values of $r$. Hence, we only have to make $\varepsilon$ small enough and $r$ big enough, such that $\bP$ becomes positive definite.

ii) $\Rightarrow$ iii):
		From ii) it is clear that we have a special solution for $n = d$. Let $\widetilde{\bA}^{(d)}, \bB^{(d)}, \bC^{(d)}, \bP^{(d)}$ be this special solution. This solution can be extended to a solution for arbitrary dimension $n \geq d$ by setting
	\begin{align*}
		\widetilde{\bA} &= \begin{pmatrix}
		\widetilde{\bA}^{(d)} & 0_{d\times n-d}\\
		0_{n-d\times d}& 0_{n-d\times n-d}
		\end{pmatrix},
		\bB = \begin{pmatrix}
		\bB^{(d)}\\
		0_{n-d\times d}
	\end{pmatrix},\\
	\bC &= \begin{pmatrix}
		\bC^{(d)} & 0_{d \times n - d}.
	\end{pmatrix},
	\bP_{22} = \bP_{22}^{(d)},\\
	\bP_{11} &= \begin{pmatrix}
		\bP_{11}^{(d)} & 0_{d\times n-d}\\
		0_{n-d\times d}& \bI_{n-d}
	\end{pmatrix}
	\bP_{12} = \begin{pmatrix}
		\bP_{12}^{(d)}\\
		0_{n-d\times d}
	\end{pmatrix}.
	\end{align*}
	Showing that these values satisfy the constraints and the LMI of Theorem \ref{thm:4.5} is straight forward.

iii) $\Rightarrow$ iv):
		As stated in Theorem \ref{thm:4.6SynIneq}, the constraints and matrix inequality of this theorem are equivalent to the matrix inequality of Theorem \ref{thm:4.5} in the case $n \geq 3d$.

 iv) $\Rightarrow$ i):
	If Theorem \ref{thm:4.6SynIneq} admits a solution, then there exists an optimizer which satisfies the conditions of Theorem \ref{thm:4.5} and thus a solution to Problem \ref{problem2}. By Theorem \ref{thm:ConvTrInv} this solution would also solve Problem \ref{problem1}, which can only be solved, if any function $f \in S(\bM,\bL)$ has a fixed point. If any function $f \in S(\bM,\bL)$ has a fixed point, then holds $\bM \preceq_c \bL$ and $\bM$ and $\bL$ are non-singular by Theorem \ref{thm:ExMinCoCa}. \qed

\subsection{Proof of Lemma  \ref{lem:FDI-LMI-connection} }

First, notice that $\sigma (\widetilde{\bA}) \subsEq \bbC_{|z| < \rho}$ is implied by the Lyapunov inequality $\widetilde{\bA}^T\bP_{11} \widetilde{\bA} - \rho^2 \bP_{11} \prec 0$, which is the left upper block of the matrix inequality \eqref{eq:4.5_1}.
In this proof, we show the equivalence of the FDI \eqref{eq:FDI} of Theorem \ref{thm:3.1ExpStabIQC} and the matrix inequality \eqref{eq:4.5_1}. Therefore, notice, the multiplier $\bPi$ from \eqref{eq:3.1IQC_1} can be factorized into $\bPi(z) = \bPsi (z)^* \bR \bPsi(z)$ with
\begin{align*}
    \resizebox{\linewidth}{!}{$
    \bPsi(z) = \begin{pmatrix}
	(l-m)(1-\lambda z^{-1}) & z^{-1}\lambda \bI_d\\
	0 & \bI_d
	\end{pmatrix},
    \bR = \begin{pmatrix}
    0 & \bI_d\\
    \bI_d & -2 \bI_d
    \end{pmatrix}.$}
\end{align*}
The goal is to apply the discrete-time KYP-Lemma to \eqref{eq:FDI}. Therefore, a realization of the following concatenation of  $\bG$ and $\bPsi$ is needed:
{\small
\begin{align*}
	\bPsi(z)\begin{pmatrix}	
	\bG(z)\\
	\bI_d
	\end{pmatrix}
	&=
	\begin{pmatrix}
	(l-m)(1-\lambda z^{-1}) & z^{-1}\lambda \bI_d\\
	0 & \bI_d
	\end{pmatrix}
	\begin{pmatrix}	
	\bG(z)\\
	\bI_d
	\end{pmatrix}\\
	&=
	\begin{pmatrix}
	z^{-1}\lambda\bI_d + (l-m)(1-\lambda z^{-1})\bG(z)\\
	\bI_d
	\end{pmatrix}\\
	&= 
	\underbrace{\begin{pmatrix}
		\lambda \bI_d + (l-m)(z- \lambda)\bG(z) & 0\\
		0 & \bI_d
		\end{pmatrix}}_{\bH_1}
	\underbrace{\begin{pmatrix}
		z^{-1}\bI_d\\
		\bI_d
		\end{pmatrix}}_{\bH_2}
\end{align*}}
Here, $\bH_1$ is realizable, because $(z-\lambda)\bG(z)$ is realizable, because $\bG$ has a relative degree of at least one. A realization of $\bH_1$ is
\begin{align*}
	\left[
	\begin{array}{c|cc}
	\widetilde{\bA} & \bB & 0 \\
	\hline
	(l-m)(\bC \widetilde{\bA} - \lambda\bC) & (l-m)\bC \bB + \lambda\bI_d & 0\\
	0 & 0 & \bI_d
	\end{array}
	\right]
\end{align*}
and a realization of $\bH_2$ is
\begin{align*}
	\left[
	\begin{array}{c|c}
	0 & \bI_d \\
	\hline
	\bI_d & 0\\
	0 & \bI_d
	\end{array}
	\right].
\end{align*}
With these realizations, the following is a realization of the chaining of $\bH_1$ and $\bH_2$:
\begin{align*}
	%\bH_1\circ\bH_2 =
	\left[
	\begin{array}{cc|c}
	\widetilde{\bA} & \bB & 0 \\
	0 & 0 & \bI_d\\
	\hline
	(l-m)(\bC \widetilde{\bA} - \lambda\bC) & (l-m)\bC \bB + \lambda\bI_d & 0\\
	0 & 0 & \bI_d
	\end{array}
	\right].
\end{align*}
This is now also a realization of $\bPsi \begin{pmatrix}
	\bI_d\\ \bG
	\end{pmatrix}$. Satisfaction of the FDI
\begin{align*}
	\begin{pmatrix}
	\bI_d\\
	\bG(z)
	\end{pmatrix}^* \bPsi^*(z) \bR \bPsi(z) \begin{pmatrix}
	\bI_d\\
	\bG(z)
	\end{pmatrix} &\prec 0 &\forall z \in \bbC_{|z| = \rho}
\end{align*}
is by the Kalman Yakubovic Popov Lemma (Corollary 13 of \cite{boczar2017exponential}) equivalent to existence of $\bP = \bP^T$ such that
{\small
\begin{align*}
	&\begin{pmatrix}
	\bI_n & 0 & 0\\
	0 & \bI_d & 0\\
	\widetilde{\bA} & \bB & 0\\
	0 & 0 & \bI_d
	\end{pmatrix}^T
	\begin{pmatrix}
	-\rho^2\bP_{11} & -\rho^2\bP_{12} & 0 & 0\\
	-\rho^2\bP_{21} & -\rho^2\bP_{22} & 0 & 0\\
	0 & 0 & \bP_{11} & \bP_{12}\\
	0 & 0 & \bP_{21} & \bP_{22}
	\end{pmatrix}
	(\star)\\
% 	\begin{pmatrix}
% 	\bI_n & 0 & 0\\
% 	0 & \bI_d & 0\\
% 	\widetilde{\bA} & \bB & 0\\
% 	0 & 0 & \bI_d
% 	\end{pmatrix}\\
	&+
	\begin{pmatrix}
	(l-m)(\bC \widetilde{\bA} - \lambda\bC) & (l-m)\bC \bB + \lambda\bI_d & 0\\
	0 & 0 & \bI_d
	\end{pmatrix}^T
	\bR\\
	&\begin{pmatrix}
	(l-m)(\bC \widetilde{\bA} - \lambda\bC) & (l-m)\bC \bB + \lambda\bI_d & 0\\
	0 & 0 & \bI_d
	\end{pmatrix}
\end{align*}
}
is negative definite. A quick reformulation of the above terms reveals that they are nothing but inequality \eqref{eq:4.5_1}. It can be checked that inequality \eqref{eq:4.5_1} can only have positive definite solutions $\bP$.\qed

\subsection{Proof of Lemma \ref{lem:5.1}}

We have to check whether inequality \eqref{eq:4SCLcharacterisation} holds for the Lagrangian function $L \in C^1$. 
Let therefore arbitrary values $z_1,z_2 \in \bbR^d$ and $\lambda_1,\lambda_2\in \bbR^{d_2}$ be given. The lower bound  in inequality \eqref{eq:4SCLcharacterisation} follows from
\begin{align*}
	&\frac{1}{2}
	\begin{pmatrix}
		z_1 - z_2\\
		\lambda_1 - \lambda_2
	\end{pmatrix}^T
	\begin{pmatrix}
		\bM & \bA_\mathrm{eq}^T\\
		\bA_\mathrm{eq} & 0
	\end{pmatrix}
	\begin{pmatrix}
		z_1 - z_2\\
		\lambda_1 - \lambda_2
	\end{pmatrix}\\
	&= \frac{1}{2}(z_1 - z_2)^T\bM (z_1 - z_2) + (\lambda_1 - \lambda_2)^T\bA_\mathrm{eq} (z_1 - z_2)\\
	&\leq f(z_2) - f(z_1) + (\nabla f(z_1))^T(z_1 - z_2)\\
	&\hspace{3.7mm} + (\lambda_1 - \lambda_2)^T\bA_\mathrm{eq} (z_1 - z_2)\\
	&= \underbrace{f(z_2) + \lambda_2^T(\bA_\mathrm{eq} z_2 - b_\mathrm{eq})}_{L(z_2,\lambda_2)} - \underbrace{(f(z_1) + \lambda_1^T(\bA_\mathrm{eq} z_1 - b_\mathrm{eq}))}_{L(z_1,\lambda_1)} \\ 
	&+ \underbrace{\resizebox{0.93\linewidth}{!}{$(\nabla f(z_1) + \bA_\mathrm{eq}^T\lambda_1)^T(z_1 - z_2) + (\lambda_1 - \lambda_2)^T(\bA_\mathrm{eq} z_1 - b_\mathrm{eq})$}}_{
	(\nabla L(z_1,\lambda_1))^T
	\begin{pmatrix}
		z_1 - z_2\\
		\lambda_1 - \lambda_2
	\end{pmatrix}
	}.
\end{align*}
The upper bound can be shown analogously. \qed

\subsection{Auxiliary results}
\begin{lemma}
    \label{lem:helpInequality}
    Let $\bM \preceq_c \bL$ be non-singular, symmetric matrices. Then for any $\rho > \rho \left( (\bM + \bL)^{-1} (\bL - \bM)\right)$ there exists an $r_0\in \bbR_{>0}$ such that for all real numbers $r \geq r_0$
	\begin{align*}
	\resizebox{\linewidth}{!}{$
	\bL - \bM \prec \rho^2 (\bL + \bM)\left((\bL - \bM)^\dagger + r\bPi_{\ker \bL - \bM}\right)(\bL + \bM).$}
	\end{align*}
\end{lemma}
\begin{proof}
	Let $\rho > \rho \left( (\bM + \bL)^{-1} (\bL - \bM)\right)$ be given. Define $\bPi := \bPi_{\ker \bL - \bM}$ and
	\begin{align*}
	\tilde{\rho} &:= 
	\rho \left( (\bL + \bM)^{-1}\left(\bL - \bM + \frac{1}{r}\bPi\right) \right)\\
	&\overset{(\star)}{=} \rho \left( \sqrt{\bL - \bM + \frac{1}{r}\bPi}(\bL + \bM)^{-1}\sqrt{\bL - \bM + \frac{1}{r}\bPi} \right)\\
	&=
	\left\|\sqrt{\bL - \bM + \frac{1}{r}\bPi}(\bL + \bM)^{-1}\sqrt{\bL - \bM + \frac{1}{r}\bPi}\right\|.
	\end{align*}
	Here, the equality $(\star)$ holds by a similarity transform with $\sqrt{\bL - \bM + \frac{1}{r}\bPi}$.
	This definition of $\tilde{\rho}$ implies the matrix inequality
	\begin{align*}
	\resizebox{\linewidth}{!}{$
	\tilde{\rho}^2\bI_d \succeq \left(\sqrt{\bL - \bM + \frac{1}{r}\bPi}(\bL + \bM)^{-1}\sqrt{\bL - \bM + \frac{1}{r}\bPi}\right)^2.$}
	\end{align*}
	A congruence transform with $\left(\bL - \bM + \frac{1}{r}\bPi\right)^{-\frac{1}{2}}(\bL + \bM)$ yields
	\begin{align*}
	\resizebox{\linewidth}{!}{$
	\tilde{\rho}^2 (\bL + \bM)\left((\bL - \bM)^\dagger + r\bPi\right)(\bL + \bM) \succeq \bL - \bM + \frac{1}{r}\bPi,$}
	\end{align*}
	since $\left(\bL - \bM + \frac{1}{r}\bPi\right)^{-\frac{1}{2}} \overset{\eqref{eq:InverseWithPseudoInverse2}}{=} \sqrt{(\bL - \bM)^\dagger + r\bPi}$.
	By the expression of $\tilde{\rho}$ through the spectral norm and the continuity of the norm, it is clear that $\tilde{\rho}$ converges to $\rho \left( (\bM + \bL)^{-1} (\bL - \bM)\right)$ for $r \to \infty$. Hence, we can choose $r$ large enough, such that $\tilde{\rho}$ is small than $\rho$ and thus,
	\begin{align*}
	\bL - \bM & \preceq
	\resizebox{0.80\linewidth}{!}{$
	\tilde{\rho}^2 (\bL + \bM)\left((\bL - \bM)^\dagger + r\bPi_{\ker\bL - \bM}\right)(\bL + \bM)$}\\
	& \prec 
	\resizebox{0.80\linewidth}{!}{$
	\rho^2 (\bL + \bM)\left((\bL - \bM)^\dagger + r\bPi_{\ker\bL - \bM}\right)(\bL + \bM)$}.
	\end{align*}
	Since increasing $r$ corresponds to adding a positive definite term to the right hand side of this inequality, the inequality remains valid for larger values of $r$.
\end{proof}

%\ce{Check if result can be found in a book.}
\begin{lemma}[Congruence Lemma]
    \label{lem:congruence}
    Let $\bM,\bL \in \bbR^{d\times d}$ be two symmetric matrices such that there exists a positive definite matrix $\bP = \bP^T \in \bbR^{d\times d}$ with
    \begin{align*}
	    \bM \bP \bL + \bL\bP\bM \succ 0.
    \end{align*}
    Then $\bM$ and $\bL$ are congruent, i.e. there exists a non-singular matrix $\bT$ such that $\bT^T\bM\bT = \bL$.
\end{lemma}
\begin{proof}
	By $\bP$ being positive definite, there exists a symmetric positive definite matrix $\sqrt{\bP} \in \bbR^{d\times d}$ with $\sqrt{\bP}^2 = \bP$. A congruence transform with $\sqrt{\bP}$ yields
	\begin{align}
		\sqrt{\bP}\bM \sqrt{\bP}\sqrt{\bP} \bL\sqrt{\bP} + \sqrt{\bP}\bL\sqrt{\bP}\sqrt{\bP}\bM\sqrt{\bP} \succ 0. \label{eq:congruence_1}
	\end{align}
	The matrices $\widetilde{\bM}:=\sqrt{\bP}\bM\sqrt{\bP}$ and $\widetilde{\bL} := \sqrt{\bP}\bL\sqrt{\bP}$ are congruent to $\bM$ and $\bL$. Hence, it is sufficient to show that the matrices $\widetilde{\bM}$ and $\widetilde{\bL}$ are congruent.\\
	Therefore, let $\bT$ be an orthogonal matrix, such that
	\begin{align*}
	\bT^T\widetilde{\bM}\bT = \begin{pmatrix}
	\bD_1 & 0\\
	0 & \bD_2
	\end{pmatrix},
	\end{align*}
	where $\bD_1$ is the diagonal matrix of all positive eigenvalues of $\widetilde{\bM}$ and $\bD_2$ is the matrix of all negative eigenvalues of $\widetilde{\bM}$. Now, a congruence transform with $\bT$ can be applied to \eqref{eq:congruence_1}:
	{\small
	\begin{align*}
	0&\prec \bT^T\widetilde{\bL}\widetilde{\bM}\bT + \bT^T \widetilde{\bM} \widetilde{\bL} \bT\\
	&=\underbrace{\bT^T\widetilde{\bL} \bT}_{ :=\bE^T} \bT^T\widetilde{\bM}\bT + \bT^T \widetilde{\bM} \bT \underbrace{\bT^T \widetilde{\bL} \bT}_{:=\bE}\\
	&=
	\begin{pmatrix}
	\bE_{11} & \bE_{12}\\
	\bE_{21} & \bE_{22}
	\end{pmatrix}
	\begin{pmatrix}
	\bD_1 & 0\\
	0 & \bD_2
	\end{pmatrix}
	+
	\begin{pmatrix}
	\bD_1 & 0\\
	0 & \bD_2
	\end{pmatrix}
	\begin{pmatrix}
	\bE_{11} & \bE_{12}\\
	\bE_{21} & \bE_{22}
	\end{pmatrix}.
	\end{align*}}
	From this inequality, one can read off
	\begin{align*}
	\bE_{11}\bD_1 + \bD_1\bE_{11} \succ 0, \hspace{10mm} \bE_{22}\bD_2 + \bD_2\bE_{22} \succ 0
	\end{align*}
	from the diagonal blocks. Hence, by the Lyapunov inequality, $\bE_{11} \succ 0$ and $\bE_{22} \prec 0$.
	Now, $\bE$ is positive definite on the subspace corresponding to $\bE_{11}$ and negative definite on the subspace corresponding to $\bE_{22}$. Consequently, $\bE$ has exactly $\dim \bE_{11} = \dim \bD_1$ positive and exactly $\dim \bE_{22} = \dim \bD_2$ negative eigenvalues according to Sylvester's law of inertia. Thus the matrices $\bM$ and $\bL$, which are congruent to $\bD$ and $\bE$, are congruent to each other.
\end{proof}

\end{document}